\documentclass[10pt,letterpaper]{amsart}
\usepackage{charter}
\usepackage[OT2,OT1]{fontenc}
\usepackage[latin9]{luainputenc}
\usepackage{array}
\usepackage{float}
\usepackage{bm}
\usepackage{multirow}
\usepackage{amstext}
\usepackage{amsthm}
\usepackage{amssymb}
\usepackage{cancel}
\usepackage{stmaryrd}

\makeatletter

\@ifundefined{pageheight}{\let\pageheight\pdfpageheight}{}
\@ifundefined{pagewidth}{\let\pagewidth\pdfpagewidth}{}
\pageheight\paperheight
\pagewidth\paperwidth

\providecommand{\tabularnewline}{\\}

\numberwithin{equation}{section}
\numberwithin{figure}{section}
\theoremstyle{plain}
\newtheorem{thm}{\protect\theoremname}[section]
\theoremstyle{definition}
\newtheorem{defn}[thm]{\protect\definitionname}
\theoremstyle{plain}
\newtheorem{lem}[thm]{\protect\lemmaname}
\theoremstyle{remark}
\newtheorem{rem}[thm]{\protect\remarkname}
\theoremstyle{plain}
\newtheorem{prop}[thm]{\protect\propositionname}
\theoremstyle{plain}
\newtheorem{cor}[thm]{\protect\corollaryname}

\usepackage{amsfonts}
\usepackage{mathrsfs}
\usepackage{graphicx}
\usepackage{amscd}

\newcommand{\cyr}{
\renewcommand\rmdefault{wncyr} \renewcommand\sfdefault{wncyss} \renewcommand\encodingdefault{OT2} \normalfont
\selectfont
}
\DeclareTextFontCommand{\textcyr}{\cyr}    

   \def\@settitle
{
\begin{center} \baselineskip14\p@\relax \LARGE\@title \end{center}
}

\usepackage[normalem]{ulem}
\usepackage[mathscr]{eucal}
\numberwithin{equation}{section}

\usepackage{verbatim}
\usepackage{amscd}
\usepackage{bm}\usepackage[mathscr]{eucal}
\usepackage[all]{xy}
\usepackage{enumerate}
\usepackage{color}
\usepackage{colortbl}

\input xyimport.tex

\title{Duality Related with Key Varieties of $\mQ$-Fano 3-folds. I}

\author{Hiromichi Takagi}

\address{Department of Mathematics, Gakushuin University, 
Mejiro, Toshima-ku, Tokyo 171-8588, Japan}
\email{hiromici@math.gakushuin.ac.jp}

\newcommand{\sA}{\mathcal{A}}
\newcommand{\sB}{\mathcal{B}}

\newcommand{\sE}{\mathcal{E}}

\newcommand{\sH}{\mathcal{H}}

\newcommand{\sO}{\mathcal{O}}
\newcommand{\sP}{\mathcal{P}}

\newcommand{\sQ}{\mathcal{Q}}

\newcommand{\sU}{\mathcal{U}}

\newcommand{\mC}{\mathbb{C}}

\newcommand{\rG}{\mathrm{G}}

\newcommand{\mP}{\mathbb{P}}
\newcommand{\mQ}{\mathbb{Q}}

\newcommand{\Pic}{\mathrm{Pic}\,}

\newcommand{\Sing}{\mathrm{Sing}\,}

\newcommand{\rank}{\mathrm{rank}\,}

\numberwithin{equation}{section}

 \newcounter{myparagraph}[subsection]

\makeatother

\providecommand{\corollaryname}{Corollary}
\providecommand{\definitionname}{Definition}
\providecommand{\lemmaname}{Lemma}
\providecommand{\propositionname}{Proposition}
\providecommand{\remarkname}{Remark}
\providecommand{\theoremname}{Theorem}

\begin{document}
\maketitle 
\begin{abstract}
In our previous paper \cite{Tak2}, we show that any prime $\mQ$-Fano
3-folds $X$ with only $1/2(1,1,1)$-singularities in certain 5 classes
can be embedded as linear sections into bigger dimensional $\mQ$-Fano
varieties called key varieties, where each of the key varieties is
constructed from certain data of the Sarkisov link staring from the
blow-up at one $1/2(1,1,1)$-singularity of $X$. In this paper, we
introduce varieties associated with the key varieties which are dual
in a certain sense. As an application, we interpret a fundamental
part of the Sarkisov link for each $X$ as a linear section of the
dual variety. In a natural context describing the variety dual to
the key variety of $X$ of genus 5 with one $1/2(1,1,1)$-singularity,
we also characterize a general canonical curve of genus 9 with a $g_{7}^{2}$. 
\end{abstract}

\maketitle
\markboth{Duality}{Hiromichi Takagi} {\small{}{\tableofcontents{}}}{\small\par}

2020\textit{ Mathematics subject classification}: 14J45, 14E30, 14H45.


\section{\textbf{Introduction}}

\subsection{Background}

This is a companion paper to \cite{Tak2}. 

In this paper, we work over $\mC$, the complex number field. For
a vector bundle $\sE$ on a variety $X,$ the notation $\mP_{X}(\sE)$
or simply $\mP(\sE)$ is just the projectivization (\textit{We don't
use the Grothendieck notation}). 

A projective variety $X$ is called a \textit{$\mQ$-Fano variety}
if $X$ has only terminal singularities and $-K_{X}$ is ample. A
$\mQ$-Fano variety $X$ is called \textit{prime} if $-K_{X}$ generates
the group of numerical equivalence classes of $\mQ$-Cartier divisors
on $X$. 

In \cite{Tak2}, we study prime $\mQ$-Fano 3-folds $X$ in the 5
classes No.1.1,~1.4,~1.9,~1.10, and 1.13 among \cite[Table 1]{Tak1},
and construct key varieties for them (see Theorem \ref{thm:IIImain}
below for the precise statement). 

\subsection{Duality for the key varieties \label{subsec:Duality}}

Let $X$ be a smooth prime Fano 3-fold of genus 9. Fano \cite[p.207-208]{Fa}
and Iskovskih \cite{Is} showed that the double projection of $X$
from a line ends with the blow-up of $\mP^{3}$ along a non-hyperelliptic
smooth curve $C$ of genus 3 and degree 7. Mukai \cite{Mu1,Mu6} showed
that $X$ is a linear section of the symplectic Grassmanian ${\rm Sp}(3,6)$.
Note that the projectively dual variety of ${\rm Sp(3,6)}$ is a quartic
hypersurface $\sH$. He also showed that the canonical model of $C$
is a linear section of $\sH$ (\cite{Mu5}). He obtained similar results
in the case of genus $7$ or $10.$ Hence Mukai revealed that the
projective duality amplifies the geometry of the Sarkisov link of
smooth Fano 3-folds. The main result of this paper concerns suitable
dual varieties to our key varieties and is modeled on these results
on duality by Mukai.

\subsection{Prime $\mQ$-Fano 3-fold and Sarkisov link}

In this subsection, we quickly review the result of \cite{Tak2} while
introducing notation which is needed in this paper. The data of $\mQ$-Fano
3-folds $X$ in the 4 classes 1.4,~1.9,~1.10, and 1.13 are summarized
in the following table:

\begin{table}[H]
\begin{tabular}{|c|c|c|c|c|c|}
\hline 
Name & No. & $g(X)$ & $\deg C$ & $g(C)$ & $X'$\tabularnewline
\hline 
\hline 
genus 5 & 1.4 & $5$ & $9$ & $9$ & $\mP^{3}$\tabularnewline
\hline 
genus 6, $\mathtt{C}$-type & 1.9 & $6$ & $3$ & $0$ & $B_{3}$\tabularnewline
\hline 
genus 6, $\mathtt{Q}$-type & 1.10 & $6$ & $9$ & $6$ & $Q^{3}$\tabularnewline
\hline 
genus 8 & 1.13 & $8$ & $7$ & $2$ & $B_{5}$\tabularnewline
\hline 
\end{tabular}
\end{table}
\noindent In the first column of the table, we rename the 4 classes.
The number $g(X)$ in the third column of the table is the \textit{genus}
of $X$ defined to be $h^{0}(-K_{X})-2.$ We explain the data in 4th--6th
column below. We recall that each $X$ in the 4 classes has only one
$1/2(1,1,1)$-singularity. We classify them in \cite{Tak1} by constructing
the following Sarkisov links:

\begin{equation}\label{eq:Sarkisov} \xymatrix{& Y\ar@{-->}[r]\ar[dl]_f & Y'\ar[dr]^{f'}\\
X & & & X',}
\end{equation}where $f\colon Y\to X$ is the blow-up of $X$ at the unique $1/2(1,1,1)$-singularity,
$Y\dashrightarrow Y'$ is a flop, and $f'$ is the blow-up of a smooth
$\mQ$-Fano $3$-fold $X'$ along a smooth curve $C$ with the genus
$g(C)$ and the degree $\deg C$ as in the 4th and 5th column of the
table, where the degree of $C$ is measured by the primitive Cartier
divisor on $X'$. In the 6th column, $B_{3}$ is a smooth cubic 3-fold
in $\mP^{4}$, $B_{5}$ is a codimension $3$ smooth linear section
of ${\rm G}(2,5)$, and $Q^{3}$ is a smooth quadric 3-fold. 

For a prime $\mQ$-Fano 3-fold $X$ of No.1.1, we rename it a prime
$\mQ$-Fano 3-fold of genus 4. Note that $X$ has two $1/2(1,1,1)$-singularities.
For such an $X$, we construct the following diagram in \cite{Tak2}:\begin{equation}\label{eq:SarkisovG4} \xymatrix{& Z\ar@{-->}[r]\ar[dl]_g & Z'\ar[dr]^{g'}\\
X & & & B_6,}
\end{equation}where $g\colon Z\to X$ is the blow-up of $X$ at the two $1/2(1,1,1)$-singularities,
$Z\dashrightarrow Z'$ is a flop, and $g'$ is the blow-up of $B_{6}:=\mP(\Omega_{\mP^{2}}^{1})$
along a smooth curve $C$ with the genus $8$ and the degree $14$.

In this paper, we call collectively the diagram (\ref{eq:SarkisovG4})
in the genus 4 case and the diagram (\ref{eq:Sarkisov}) in the other
cases \textit{the basic diagram} (we also keep the name the Sarkisov
link for the diagram (\ref{eq:Sarkisov})).

We say that a projective variety $X$ is a \textit{linear section}
of a projective variety $\Sigma$ with respect to a linear system
$|M_{\Sigma}|$ if it holds that $X=\Sigma\cap D_{1}\cap\dots\cap D_{k}$
for $k=\dim\Sigma-\dim X$ and some $D_{1},\dots,D_{k}\in|M_{\Sigma}|$.
We usually do not mention the linear system $|M_{\Sigma}|$ if $M_{\Sigma}$
generates the group of the numerical equivalence classes of $\mQ$-Cartier
divisors on $\Sigma$. We can say that the main result of \cite{Tak2}
as follows is a classification of $\mQ$-Fano 3-folds in the 5 classes
in different nature to that in \cite{Tak1}. 
\begin{thm}[Embedding theorem \cite{Tak2}]
 \label{thm:IIImain} For each one of the $5$ classes, there is
a unique rational $\mQ$-Fano variety $\Sigma$ of Picard number $1$
such that any prime $\mQ$-Fano $3$-fold $X$ in the class is a linear
section of $\Sigma$. The $\mQ$-Fano varieties $\Sigma$ are of $11$-,
$12$-, $9$-, $8$-, and $5$-dimensional for $X$ of genus $4$,
$5$, of genus $6$ and $\mathtt{Q}$-type, of genus $6$ and $\text{\ensuremath{\mathtt{C}}}$-type,
and of genus $8$, respectively. 
\end{thm}

For a prime $\mQ$-Fano 3-fold $X$ in each of the 5 classes, we will
call the variety $\Sigma$ \textit{the key variety} for $X$. 

\subsection{Main result}

Through the constructions of the key varieties $\Sigma$, we obtain
their birational models which are projective bundles over Fano manifolds
$S$ associated with certain vector bundles $\sE$ such that $\sE^{*}$
is globally generated. In this paper, the projective bundle $\mP(\sE)$
in each case is more important than the key variety itself. 

We set $V_{\sE}:=H^{0}(S,\sE^{*})^{*}$. Following \cite[Sect.8]{Ku4},
another vector bundle $\sE^{\perp}$ on $S$ is defined by the following
exact sequence:
\[
0\to\sE^{\perp}\to(V_{\sE})^{*}\otimes\sO_{S}\to\sE^{*}\to0.
\]
Here is the table of the data as we have mentioned with notation and
conventions below:

\begin{table}[H]
Table 1 %
\begin{tabular}{|c|c|c|}
\hline 
\multirow{3}{*}{$g=4$} & $S$ & $B_{6}\subset\mP(\mathrm{S}^{-1,0,1}U^{3})\simeq\mP^{7}$\tabularnewline
\cline{2-3} \cline{3-3} 
 & $\sE$ & $p_{1}^{*}\sO_{\mP((U^{3})^{*})}(-1)\oplus p_{2}^{*}\sO_{\mP(U^{3})}(-1)\oplus\Omega_{\mP(\mathrm{S}^{-1,0,1}U^{3})}^{1}(1)|_{B_{6}}$\tabularnewline
\cline{2-3} \cline{3-3} 
 & $\sE^{\perp}$ & $p_{1}^{*}\Omega_{\mP((U^{3})^{*})}^{1}(1)\oplus p_{2}^{*}\Omega_{\mP(U^{3})}^{1}(1)\oplus\sO_{\mP(\mathrm{S}^{-1,0,1}U^{3})}(-1)|_{B_{6}}$\tabularnewline
\hline 
\multirow{3}{*}{$g=5$} & $S$ & $\mP(U^{4})\simeq\mP^{3}$\tabularnewline
\cline{2-3} \cline{3-3} 
 & $\sE$ & $U^{3}\otimes\Omega_{\mP(U^{4})}^{1}(1)\oplus\sO_{\mP(U^{4})}(-1)$\tabularnewline
\cline{2-3} \cline{3-3} 
 & $\sE^{\perp}$ & $(U^{3})^{*}\otimes\sO_{\mP(U^{4})}(-1)\oplus\Omega_{\mP(U^{4})}^{1}(1)$\tabularnewline
\hline 
\multirow{3}{*}{$g=6$, $\mathsf{Q}$-type} & $S$ & $Q^{3}\subset\mP(U^{5})$\tabularnewline
\cline{2-3} \cline{3-3} 
 & $\sE$ & $\sU|_{Q^{3}}\oplus\sO_{Q^{3}}(-1)\oplus\Omega_{\mP(U^{5})}^{1}(1)|_{Q^{3}}$\tabularnewline
\cline{2-3} \cline{3-3} 
 & $\sE^{\perp}$ & $\sQ^{*}|_{Q^{3}}\oplus\Omega_{\mP(U^{5})}^{1}(1)|_{Q^{3}}\oplus\sO_{Q^{3}}(-1)$\tabularnewline
\hline 
\multirow{3}{*}{$g=6$, $\mathsf{C}$-type} & $S$ & $\widehat{A}_{\mathtt{C}}$\tabularnewline
\cline{2-3} \cline{3-3} 
 & $\sE$ & $a^{*}\sO_{A_{\mathtt{C}}}(-1)\oplus b^{*}\Omega_{\mP(U^{5})}^{1}(1)$\tabularnewline
\cline{2-3} \cline{3-3} 
 & $\sE^{\perp}$ & $a^{*}(\Omega_{\mP(U^{8})}^{1}(1)|_{A_{\mathtt{C}}})\oplus b^{*}\sO_{\mP(U^{5})}(-1)$\tabularnewline
\hline 
\multirow{3}{*}{$g=8$} & $S$ & $B_{5}\subset\mP(U^{7})$\tabularnewline
\cline{2-3} \cline{3-3} 
 & $\sE$ & $\sU|_{B_{5}}\oplus\sO_{B_{5}}(-1)$\tabularnewline
\cline{2-3} \cline{3-3} 
 & $\sE^{\perp}$ & $\sQ^{*}|_{B_{5}}\oplus\Omega_{\mP(U^{7})}^{1}(1)|_{B_{5}}$\tabularnewline
\hline 
\end{tabular}
\end{table}

\begin{itemize}
\item $U^{i}$: a $i$-dimensional vector space.
\item (For the genus 4 case) We consider 
\[
B_{6}=\{{\empty}^{t}\!\bm{y}\bm{x}=0\}\subset\mP((U^{3})^{*})\times\mP(U^{3}),
\]
where $\text{\ensuremath{{\empty}^{t}\!\bm{y}}\ensuremath{\ensuremath{\in}(\ensuremath{U^{3})^{*}}}}$
and $\bm{x}\in U^{3}$ are considered as row and column vectors respectively.
We also identify $B_{6}$ with its image by the Segre embedding 
\begin{align*}
S\colon\mP((U^{3})^{*})\times\mP(U^{3}) & \hookrightarrow\mP((U^{3})^{*}\otimes U^{3})\\{}
[\bm{y}]\times[\bm{x}] & \mapsto[\bm{y}\otimes\bm{x}].
\end{align*}
Then $B_{6}$ spans $\mP(\mathrm{S}^{-1,0,1}U^{3})$, where $\mathrm{S}^{-1,0,1}U^{3}$
is the $8$-dimensional irreducible component of $(U^{3})^{*}\otimes U^{3}$
as ${\rm SL}(U^{3})$-representation space. We denote the natural
projections by $p_{1}\colon B_{6}\to\mP((U^{3})^{*})$ and $p_{2}\colon B_{6}\to\mP(U^{3})$,
and set $\sO_{B_{6}}(1,0):=p_{1}^{*}\sO_{\mP((U^{3})^{*})}(1)$ and
$\sO_{B_{6}}(0,1):=p_{2}^{*}\sO_{\mP(U^{3})}(1)$. 
\item $\sU:$ the universal subbundle of rank $2$ on $\rG(2,n)$, $\sQ:$
the universal quotient bundle of rank $n-2$ on $\rG(2,n)$. 
\item (For the case of genus 6 and $\mathtt{C}$-type) Let $A_{\mathtt{C}}$
be a smooth 4-dimesional linear section of ${\rm G}(2,5)$; 
\[
A_{\mathtt{C}}={\rm G}(2,5)\cap\mP(U^{8}).
\]
By \cite{Fuj}, $A_{\mathtt{C}}$ is unique up to isomorphism, and
has a unique plane $\Pi$ such that, for the blow-up $a\colon\widehat{A}_{\mathtt{C}}\to A_{\mathtt{C}}$
along $\Pi$, there exists a morphism $b\colon\widehat{A}_{\mathtt{C}}\to\mP(U^{5})$
which is the blow-up along a twisted cubic $\gamma_{\mathtt{C}}$. 
\end{itemize}
Let 
\[
\overline{\Sigma}^{*}\subset\mP((V_{\sE})^{*})
\]
 be the image of $\mP(\sE^{\perp})$ by the tautological linear system.
The main result of this paper asserts the relationship between the
basic diagram and a linear section of $\overline{\Sigma}^{*}$ in
each case:
\begin{thm}
\label{thm:main1} The following assertions hold:\vspace{3pt}

\noindent $(1)$ In the case of genus 6 and $\mathtt{C}$-type, $\overline{\Sigma}^{*}$
is a cubic $11$-fold, and the cubic 3-fold $X'$ appearing in (\ref{eq:Sarkisov})
is a linear section of $\overline{\Sigma}^{*}$ (Theorem \ref{thm:main2}
(2), Proposition \ref{prop:Dual Cubic} and Corollary \ref{cor:Cubic 3-fold}).

\vspace{3pt}

\noindent $(2)$ In the genus 8 case, $\overline{\Sigma}^{*}=\mP((V_{\sE})^{*})\simeq\mP^{11}$
and the map $\mP(\sE^{\perp})\to\overline{\Sigma}^{*}$ is a generically
finite double cover branched along a sextic hypersurface. The canonical
map of the curve $C\subset X'$ of genus $2$ can be identified with
the restriction of $\mP(\sE^{\perp})\to\overline{\Sigma}^{*}$ over
a line in $\overline{\Sigma}^{*}$ (Theorem \ref{thm:trinity}, Proposition
\ref{prop:Genus 8 dual} and Corollary \ref{cor:genus 2}). 

\vspace{3pt}

\noindent $(3)$ In each of the other cases, the canonical model
of the curve $C\subset X'$ is a linear section of $\overline{\Sigma}^{*}$
(Theorem \ref{thm:trinity} and Proposition \ref{prop:DualG4} (the
genus $4$ case), Theorem \ref{thm:trinity} and Proposition \ref{prop:G5 Dual}
(the genus $5$ case) and Theorem \ref{thm:trinity} and the explanation
as in the section \ref{sec:genus 6Q} (the case of genus $6$ and
$\mathtt{Q}$-type)). 
\end{thm}

The results of Mukai which we have mentioned in this subsection are
developed in perspective of derived category by Kuznetsov (\cite{Ku1,Ku2,Ku3,Ku4}).
Our result mentioned in this subsection can be interpreted by linear
duality \cite[Sect.8]{Ku4}, which is a special important case of
Kuznetsov's theory of homological projective duality.

\subsection{Classification of algebraic curves}

In the series of works \cite{Mu2,Mu3,Mu4,Mu7,MuId}, Mukai, partly
with Ide in the genus 8 case, has been relating generality conditions
of algebraic curves (gonality, Clifford index, Brill-Noether condition)
with key variety descriptions of them. For example, he showed in \cite{Mu3}
that a curve $C$ of genus 8 has no $g_{7}^{2}$ if and only if $C$
is a linear section of ${\rm G(2,6)}.$ As for curves of genus 8,
he, partly with Ide, completed this type of equivalence in any case
(\cite{Mu2,Mu7,MuId}). We refer to \cite{Mu2,Mu3,Mu4,Mu7} for the
results about curves of different genus. 

In this paper, we give a contribution in this direction as follows:
\begin{thm}[=Corollary \ref{cor:curveG9V2}]
 Let $C$ be a smooth curve of genus $9$. The following assertions
$({\rm a})$ and $({\rm b})$ are equivalent: 

\vspace{3pt}

\noindent $({\rm a})$ There exists a birational morphism $\iota_{1}$
from $C$ to a septic plane curve $C_{1}$ with only double points
and an isomorphism $\iota_{2}\colon C\to C_{2}$ to a space curve
$C_{2}$ of degree $9$ such that $\iota_{1}^{*}\sO_{C_{1}}(1)+\iota_{2}^{*}\sO_{C_{2}}(1)=K_{C}.$ 

\vspace{3pt}

\noindent $({\rm b})$ $C$ is isomorphic to a linear section of
$\overline{\Sigma}^{*}$ associated with prime $\mQ$-Fano $3$-fold
of genus $5$.
\end{thm}

We note that a curve of genus 9 with condition (a) has Clifford index
3 and admit a $g_{7}^{2}$ (cf.~\cite{Sa}) but the converse is not
true in general. We refer Remark \ref{rem:G9G27} to more detailed
explanations as for this. 

We also reproduce a result of Mukai in \cite{Mu7} about a curve of
genus 8 while interpreting the key variety of the curve as the dual
to the key variety of prime $\mQ$-Fano 3-folds of genus 4 (Corollary
\ref{cor:curveG8}). 

\vspace{5pt}

\noindent\textbf{ Notation and Conventions} 
\begin{itemize}
\item \textit{Tautological line bundle}: Setting $\mathcal{P}=\mP(\sE)$,
we often denote by $\sO_{\sP}(1)$, or $H_{\sP}$ the tautological
line bundle associated to the vector bundle $\sE$. 
\item \textit{Point of a projective space}: Let $V$ be a vector space.
For a nonzero vector $\bm{x}\in V$ and a 1-dimensional subspace $V^{1}\subset V$,
we denote by $[\bm{x}]$ and $[V^{1}]$ the point of $\mP(V)$ corresponding
to $\bm{x}$ and $V^{1}$ respectively.
\item \textit{Cartier divisor and invertible sheaf}: We sometimes abuse
notation of a Cartier divisor and an invertible sheaf. For example,
we sometimes use the expression like $D=f^{*}\sO_{X}(1)$.
\end{itemize}
\textbf{\noindent Acknowledgment}: I am grateful to Professor Shinobu
Hosono for his encouragement while writing this paper. Through previous
collaborations with him, my understanding about linear duality were
deepened and this led me to Theorem \ref{thm:main2}. A strong motivation
to obtain results in this paper came from Professor Mukai's another
side of works \cite{Mu3,Mu4}. From personal conversations with him,
I learned a lot of things about this. I appreciate him giving me a
lot of ideas generously. This work is supported in part by Grant-in
Aid for Scientific Research (C) 16K05090.

\section{\textbf{Duality related with key varieties \label{sec:LinearDuality}}}

\subsection{Generalities on vector bundle}

We follow \cite[Sect.8]{Ku4} but we only consider the situation as
in the subsection \ref{subsec:Duality}. 

We denote by $\pi\colon\mP(\sE)\to S$ the natural projection, and
by $\varphi\colon\mP(\sE)\to\mP(V_{\sE})$ the morphism defined by
the tautological linear system of $\mP(\sE)$. Similarly, we denote
by $\sigma\colon\mP(\sE^{\perp})\to S$ the natural projection, and
by $\psi\colon\mP(\sE^{\perp})\to\mP((V_{\sE})^{*})$ the morphism
defined by the tautological linear system of $\mP(\sE^{\perp})$.
It should be convenient to keep these in mind as in the following
diagram:

\begin{equation}\label{eq:LinearDual} \xymatrix{& \mP(\sE)\ar[dl]_\varphi\ar[dr]^\pi & & \mP(\sE^\perp)\ar[dr]^{\psi}\ar[dl]_\sigma\\
\mP(V_\sE) & & S & & \mP((V_\sE)^*).}
\end{equation}
\begin{defn}
Let $\Lambda$ be a subspace of $(V_{\sE})^{*}$ of dimension $l$.
We set 
\[
\mP(\sE)_{\Lambda}:=\mP(\sE)\times_{\mP(V_{\sE})}\mP(\Lambda^{\perp}),\,\mP(\sE^{\perp})_{\Lambda}:=\mP(\sE^{\perp})\times_{\mP((V_{\sE})^{*})}\mP(\Lambda).
\]
We say that $\mP(\sE)_{\Lambda}$ and $\mP(\sE^{\perp})_{\Lambda}$
are \textit{mutually orthogonal linear section of $\mP(\sE)$ and
$\mP(\sE^{\perp})$} respectively if the codimension of $\mP(\sE)_{\Lambda}$
in $\mP(\sE)$ is equal to $l$ and the codimension of $\mP(\sE^{\perp})_{\Lambda}$
in $\mP(\sE^{\perp})$ is equal to $\dim V_{\sE}-l$.
\end{defn}

We refer to \cite[Lem.4.1.1]{HoTak} for a proof of the following
lemma, which is elementary but plays a crucial role in the sequel:
\begin{lem}
\label{lem:HoTak} We set $r:=\rank\sE$. Let $s\in S$ be a point.
It holds that $\dim(\sE_{s}\cap\Lambda^{\perp})=\dim(\sE_{s}^{\perp}\cap\Lambda)+r-l$.
\end{lem}

\subsection{Linear sections of $\mP(\sE)$ and $\mP(\sE^{\perp})$, and the basic
diagram}

In the following theorem, we interpret a part of the basic diagram
$(\ref{eq:Sarkisov})$ or $(\ref{eq:SarkisovG4})$ as orthogonal linear
sections of $\mP(\sE)$ and $\mP(\sE^{\perp})$.
\begin{thm}
\label{thm:main2} Let the pair $(S,\sE,\sE^{\perp})$ be as in Table
1 for each of the $5$ classes of prime $\mQ$-Fano $3$-folds. The
following assertions hold:

\vspace{3pt}

\noindent $(1)$ $\text{(On the key variety side)}$\vspace{3pt}

\noindent $(1\text{{-}}1)$ In the case of genus 4, the morphism
$Z'\to B_{6}$ appearing in the basic diagram $(\ref{eq:SarkisovG4})$
can be identified with $\pi|_{\mP(\sE)_{\Lambda}}\colon\mP(\sE)_{\Lambda}\to S$
for a linear subspace $\Lambda$ of $V^{*}$. 

\vspace{3pt}

\noindent$(1\text{{-}}2)$ In the case of genus 5, 8, or genus $6$
and $\mathtt{Q}$-type, $f'\colon Y'\to X'$ appearing in the Sarkisov
link $(\ref{eq:Sarkisov})$ can be identified with $\pi|_{\mP(\sE)_{\Lambda}}\colon\mP(\sE)_{\Lambda}\to S$
for a linear subspace $\Lambda$ of $V^{*}$. 

\vspace{3pt}

\noindent $(1\text{{-}}3)$ In the case of genus $6$ and $\mathtt{C}$-type,
$f'\colon Y'\to X'=B_{3}$ can be identified with the morphism $\mP(\sE)_{\Lambda}\to b\circ\pi(\mP(\sE)_{\Lambda})$
induced by $(b\circ\pi)|_{\mP(\sE)_{\Lambda}}$ for a linear subspace
$\Lambda$ of $V^{*}$.

\vspace{3pt}

\noindent $(2)$ $\text{(On the dual side)}$ In the case of genus
$6$ and $\mathtt{C}$-type, the morphism $f'\colon Y'\to X'=B_{3}$
can be identified with the morphism $\mP(\sE^{\perp})_{\Lambda}\to b\circ\sigma(\mP(\sE^{\perp})_{\Lambda})$
induced by $(b\circ\sigma)|_{\mP(\sE^{\perp})_{\Lambda}}$ with the
same $\Lambda$ as in $(1\text{-}3).$ In each of the other cases,
the curve $C$ appearing in the basic diagram $(\ref{eq:Sarkisov})$
or $(\ref{eq:SarkisovG4})$ is isomorphic to both $\mP(\sE^{\perp})_{\Lambda}$
and $\sigma(\mP(\sE^{\perp})_{\Lambda})$ with the same $\Lambda$
as in $(1\text{-}1)$ or $(1\text{-}2).$ 
\end{thm}

\begin{proof}
(1). In the case of genus $6$ and $\mathtt{C}$-type or of genus
8, the assertion is just a restatement of \cite[Cor. 5.18 or 3.8]{Tak2}.
In the case of genus $6$ and $\mathtt{Q}$-type or genus 4 (resp.
genus 5), the assertion follows from \cite[Cor. 5.18 and Prop. 5.22]{Tak2}
(resp. \cite[Cor. 6.16]{Tak2}) since the $\varphi_{|H_{\widehat{\Sigma}}|}$-image
of $E_{\widehat{\Sigma}}$ is disjoint from $W$ by \cite[Lem. 5.16]{Tak2}
(resp. \cite[Proof of Thm. 6.15]{Tak2}).

\vspace{3pt}

\noindent (2). \vspace{3pt}

\noindent \textbf{Cases except the case of genus 6 and $\mathtt{C}$-type}:
To treat these cases, we assume that $r-l=1$. Then, for a point $s\in S$,
we have 
\begin{equation}
\dim(\sE_{s}\cap\Lambda^{\perp})=\dim(\sE_{s}^{\perp}\cap\Lambda)+1\label{eq:dimeq}
\end{equation}
 by Lemma \ref{lem:HoTak}. 

Since $Z'\to B_{6}$ (resp. $Y'\to X'$) is the blow-up along a smooth
curve $C$ in the genus 4 case (resp. in each of the other cases),
it holds that $\dim(\sE_{s}\cap\Lambda^{\perp})=1$ (resp. $=2$)
if and only if $s\not\in C$ (resp. $s\in C$) by (1). Therefore,
by the equality (\ref{eq:dimeq}), the $\sigma$-image of $\mP(\sE^{\perp})_{\Lambda}$
is equal to $C$ and the induced morphism $\mP(\sE^{\perp})_{\Lambda}\to C$
is injective, hence is an isomorphism as desired since $C$ is smooth. 

\vspace{3pt}

\noindent \textbf{Case of genus 6 and $\mathtt{C}$-type}: To treat
this case, we assume that $r=l$. Then, for a point $s\in S=\widehat{A}_{\mathtt{C}}$,
we have 
\begin{equation}
\dim(\sE_{s}\cap\Lambda^{\perp})=\dim(\sE_{s}^{\perp}\cap\Lambda)\label{eq:EEperpequal}
\end{equation}
 by Lemma \ref{lem:HoTak}. This implies that $\sigma\left(\mP(\sE^{\perp})_{\Lambda}\right)=\pi\left(\mP(\sE)_{\Lambda}\right)$
and hence $b\circ\sigma\left(\mP(\sE^{\perp})_{\Lambda}\right)=b\circ\pi\left(\mP(\sE)_{\Lambda}\right)=X'$
by (1). Moreover, since $\mP(\sE)_{\Lambda}\to\pi\left(\mP(\sE)_{\Lambda}\right)\to b\circ\pi\left(\mP(\sE)_{\Lambda}\right)$
is the blow-up along $C$ by (1), and $b\colon\widehat{A}_{C}\to\mP(U^{5})$
is the blow-up along $C$, we have $\mP(\sE)_{\Lambda}\to\pi\left(\mP(\sE)_{\Lambda}\right)$
is an isomorphism and $\pi\left(\mP(\sE)_{\Lambda}\right)\to b\circ\pi\left(\mP(\sE)_{\Lambda}\right)$
is the blow-up along $C$. Therefore, $\mP(\sE^{\perp})_{\Lambda}\to\sigma\left(\mP(\sE^{\perp})_{\Lambda}\right)$
is an isomorphism by (\ref{eq:EEperpequal}) and $\sigma(\mP(\sE^{\perp})_{\Lambda})\to b\circ\sigma(\mP(\sE^{\perp})_{\Lambda})=X'$
is the blow-up along $C$ as desired.
\end{proof}
In the following sections, we investigate the morphism $\psi\colon\mP(\sE^{\perp})\to\mP((V_{\sE})^{*})$
and the $\psi$-image $\overline{\Sigma}^{*}$in detail in each of
the 5 cases, and show the main result Theorem \ref{thm:main1}. The
way of investigations of $\psi$ and $\overline{\Sigma}^{*}$ is similar
to that of $\mP(\sE)\to\mP(V_{\sE})$ and $\overline{\Sigma}$ as
in \cite{Tak2}. 
\begin{rem}
We can also construct the Sarkisov links related with $\psi$. For
the moment, however, we do not find an appropriate dual perspective
for them. So we do not write down them and we will revisit them in
a future.
\end{rem}

The following result is frequently used in the sequel. A proof for
this is omitted since it is elementary.
\begin{lem}
\label{lem:ABFib} Let $S$ be a projective manifold and $\sA,\sB$
vector bundles on $S$ whose dual bundles are globally generated.
Let $U_{\sA}:=H^{0}(S,\sA^{*})^{*}$and $U_{\sB}:=H^{0}(S,\sB^{*})^{*}$.
Let $p\colon\mP_{S}(\sA\oplus\sB)\to S$ be the natural morphism and
$\mu\colon\mP_{S}(\sA\oplus\sB)\to\mP(U_{\sA}\oplus U_{\sB})$ the
morphism defined by the tautological linear system $|H_{\mP(\sA\oplus\sB)}|$.
The following assertions hold:

\vspace{3pt}

\noindent $(1)$ The projective bundle $\mP_{S}(\sA\oplus\sB)$ is
contained in $\mP(U_{\sA}\oplus U_{\sB})\times S$ as a subbundle,
and the morphism $\mu$ is nothing but the composite $\mP_{S}(\sA\oplus\sB)\hookrightarrow\mP(U_{\sA}\oplus U_{\sB})\times S\to\mP(U_{\sA}\oplus U_{\sB}).$
The pull-back of $\sO_{\mP(U_{\sA}\oplus U_{\sB})}(1)$ by this morphism
is the tautological line bundle of $\mP_{S}(\sA\oplus\sB)$.

\vspace{3pt}

\noindent $(2)$ For a point $s\in S$, let $\sA_{s}$ and $\sB_{s}$
the fibers of $\sA$ and $\sB$ at $s$ respectively, which are subspaces
of $U_{\sA}$ and $U_{\sB}$ respectively. The $\mu$-image coincides
the locus
\[
\left\{ [\bm{x}+\bm{y}]\in\mP(U_{\sA}\oplus U_{\sB})\mid\exists_{s\in S},\bm{x}\in\sA_{s},\bm{y}\in\sB_{s}\right\} 
\]
 and the $\mu$-fiber over a point $[\bm{x}+\bm{y}]$ coincides with
the locus $\left\{ s\in S\mid\bm{x}\in\sA_{s},\bm{y}\in\sB_{s}\right\} .$
\end{lem}

Lemma \ref{lem:ABFib} also holds for a direct sum of three or more
vector bundles.

\section{\textbf{$\mQ$-Fano 3-fold of genus 4}}

\subsection{Descriptions of $\mP(\sE^{\perp})$ \label{subsec:DescriptionsPEperpG4}}
\begin{prop}
\label{prop:DualG4} The following assertions hold: 

\vspace{3pt}

\noindent$(1)$ $\overline{\Sigma}^{*}\subset\mP(U^{3}\oplus(U^{3})^{*}\oplus\mathrm{S}^{-1,0,1}U^{3})\simeq\mP^{12}$
is defined by the following equations $:$ 
\begin{equation}
\empty^{t}\!\bm{p}D=\empty^{t}\!\!\bm{0},D\bm{q}=\bm{0},D^{\dagger}=O,{\rm tr}D=0,\label{eq:G8DualEq}
\end{equation}
where $\empty^{t}\!\bm{p}\in U^{3}$, $\bm{q}\in(U^{3})^{*}$, and
$D\in\mathrm{S}^{-1,0,1}U^{3}$, and these are considered as a $3$-dimensional
row vector, a $3$-dimensional column vector and a traceless $3\times3$
matrix, respectively, and $D^{\dagger}$ is the adjoint matrix of
$D$.

\vspace{3pt}

\noindent $(2)$ We set $E_{\psi}:=\mP(p_{1}^{*}\Omega_{\mP((U^{3})^{*})}^{1}(1)\oplus p_{2}^{*}\Omega_{\mP(U^{3})}^{1}(1)\oplus0)$.
The morphism $\psi$ is a crepant divisorial contraction whose exceptional
locus is the divisor $E_{\psi}$.

\vspace{3pt}

\noindent $(3)$ The singular locus of $\overline{\Sigma}^{*}$ coincides
with $\mP(U^{3}\oplus(U^{3})^{*}\oplus0)$, which is the $\psi$-image
of $E_{\psi}$.

\vspace{3pt}

\noindent $(4)$ $\overline{\Sigma}^{*}$is a $7$-dimensional Fano
variety of degree $14$ with only Gorenstein canonical singularity
and with $-K_{\overline{\Sigma}^{*}}=\sO_{\overline{\Sigma}^{*}}(5)$.

\vspace{3pt}

\noindent $(5)$ The linear projection of $\mP(U^{3}\oplus(U^{3})^{*}\oplus\mathrm{S}^{-1,0,1}U^{3})$
from $\mP(U^{3}\oplus(U^{3})^{*}\oplus0)$ induces the rational map
$\overline{\Sigma}^{*}\dashrightarrow B_{6}\subset\mP(\mathrm{S}^{-1,0,1}U^{3})$.
\end{prop}

\begin{proof}
(1) Let $[W^{1}\otimes U^{1}]\in B_{6}$ be a point. The $\sigma$-fiber
over the point $[W^{1}\otimes U^{1}]$ is $\mP\left(\left((U^{3})^{*}/W^{1}\right)^{*}\oplus(U^{3}/U^{1})^{*}\oplus W^{1}\otimes U^{1}\right).$
With this description, we immediately see that $\overline{\Sigma}^{*}$
is contained in the variety defined by the equation (\ref{eq:G8DualEq}),
which we temporarily denote by $(\overline{\Sigma}^{*})'$. Let $[\empty^{t}\!\bm{p},\bm{q},D]$
be a point of $(\overline{\Sigma}^{*})'$ such that $D\not=0$. Then,
by Lemma \ref{lem:ABFib} (2), we see that the $\psi$-fiber over
$[\empty^{t}\!\bm{p},\bm{q},D]$ consists of one point $[W^{1}\otimes U^{1}]$
such that $W^{1}\otimes U^{1}=\mC D$. Therefore, $\mP(\sE^{\perp})\to(\overline{\Sigma}^{*})'$
is dominant, hence is surjective, and is also birational. 

\vspace{3pt}

\noindent (2). Since $-K_{\mP(\sE^{\perp})}=5H_{\mP(\sE^{\perp})}$,
the morphism $\psi$ is crepant. Since $\Pic\mP(\sE^{\perp})$ is
spanned by $H_{\mP(\sE^{\perp})},\sigma^{*}p_{1}^{*}\sO_{\mP((U^{3})^{*})}(1),\sigma^{*}p_{2}^{*}\sO_{\mP(U^{3})}(1)$,
any $\psi$-exceptional curve $\delta$ is positive for $\sigma^{*}p_{1}^{*}\sO_{\mP((U^{3})^{*})}(1)$
or $\sigma^{*}p_{2}^{*}\sO_{\mP(U^{3})}(1)$ since $H_{\mP(\sE^{\perp})}\cdot\delta=0$.
Since $E_{\psi}\sim H_{\mP(\sE^{\perp})}-\sigma^{*}(p_{1}^{*}\sO_{\mP((U^{3})^{*})}(1)+p_{2}^{*}\sO_{\mP(U^{3})}(1))$,
we have $E_{\psi}\cdot\delta<0$, and hence $\delta\subset E_{\psi}.$
Therefore $E_{\psi}$ is contained in the $\psi$-exceptional locus.
Since $\psi(E_{\psi})=\mP(U^{3}\oplus(U^{3})^{*})$ and $\dim E_{\psi}>\dim\mP(U^{3}\oplus(U^{3})^{*})$,
$E_{\psi}$ coincides with the $\psi$-exceptional divisor.

The assertion (3) follows from (2). As for the assertion (4), $\deg\overline{\Sigma}^{*}=14$
follows from $H_{\mP(\sE^{\perp})}^{7}=14$, the derivation of which
we omit since it is similar to the proof of Proposition \ref{prop:Dual Cubic}
(1) or \ref{prop:Genus 8 dual} below based on computations of Chern
classes of vector bundles. The remaining assertions of (4) follows
from (2) and (3). The assertion (5) immediately follows from the equation
(\ref{eq:G8DualEq}).
\end{proof}

\subsection{Curve of genus $8$ \label{subsec:Curve-of-genus8}}

Let $C$ be any smooth non-hyperelliptic curve of genus $8$ with
a non half-canonical $g_{7}^{2}$ and no $g_{4}^{1}$. By \cite{MuId},
the canonical model of $C$ is the complete intersection in $\mP^{2}\times\mP^{2}$
of three divisors of $(1,1)$-, $(2,1)$- and $(1,2)$-types. The
following corollary is just a special case of \cite[Thm.2]{Mu7} for
such a curve $C$ such that the $(1,1)$ divisor containing $C$ is
smooth. Since it is also obtained in our context naturally, we write
it down.
\begin{cor}[Curve of genus $8$]
\label{cor:curveG8} Let $C$ be a smooth curve of genus $8$. The
following are equivalent:

\vspace{3pt}
\end{cor}

\noindent $(1)$ the canonical model of $C$ is the complete intersection
in $B_{6}$ of two divisors of $(2,1)$- and $(1,2)$-types.

\vspace{3pt}

\noindent $(2)$ $C$ is projectively equivalent to a linear section
of $\overline{\Sigma}^{*}$.

\begin{proof}
The implication $(1)\Rightarrow(2)$ is a special case of \cite[Thm.2]{Mu7}.
The converse follows by reversing the discussion.
\end{proof}
\begin{rem}
The following remark should be well-known for experts: Assume that
a curve $C$ of genus 4 is the complete intersection in $\mP^{2}\times\mP^{2}$
of three divisors of $(1,1)$-, $(2,1)$- and $(1,2)$-types. Then
$C$ has a non half-canonical $g_{7}^{2}$ and no $g_{4}^{1}$. This
can be proved in a similar way to Corollary \ref{cor:nontri} and
Proposition \ref{prop:curveG9V1} $(2)\Rightarrow(1)$ below (note
that the assertion of Lemma \ref{lem:notrisec} holds also for $\mP^{2}\times\mP^{2}$
with the same proof). We add one more remark. If a curve $C$ of genus
4 is the complete intersection in $\mP^{2}\times\mP^{2}$ of three
divisors of $(1,1)$-, $(2,1)$- and $(1,2)$-type, then it has two
plane models which are the images of the first and the second projections
$\mP^{2}\times\mP^{2}\to\mP^{2}$ respectively. Two $g_{7}^{2}$ 's
which add up to $K_{C}$ are obtained as the pull-backs by the two
projections of the restrictions of lines to the two plane models.
Using the Koszul resolution of the ideal sheaf of $C\subset\mP^{2}\times\mP^{2}$,
we see that two $g_{7}^{2}$ 's are not linearly equivalent.
\end{rem}

\section{\textbf{$\mQ$-Fano 3-fold of genus 5 }}

\subsection{Descriptions of $\mP(\sE^{\perp})$ \label{subsec:DescriptionsPEperpG5}}

\begin{prop}
\label{prop:G5 Dual} The following assertions hold: 

\vspace{3pt}

\noindent$(1)$ $\overline{\Sigma}^{*}\subset\mP((U^{3})^{*}\otimes U^{4}\oplus(U^{4})^{*})\simeq\mP^{15}$
is defined by the following equations $:$ 
\begin{equation}
\empty^{t}\!\bm{p}D=0,\ \text{and}\ \rank D\leq1,\label{eq:G9DualEq}
\end{equation}
where $\empty^{t}\!\bm{p}\in(U^{4})^{*}$, and $D\in(U^{3})^{*}\otimes U^{4}$,
and these are considered as a $4$-dimensional row vector, and a $4\times3$
matrix, respectively.

\vspace{3pt}

\noindent $(2)$ Let $S_{\psi}:=\mP(0\oplus\Omega_{\mP(U^{4})}^{1}(1)).$
The morphism $\psi$ is a crepant small contraction whose exceptional
locus is $S_{\psi}$.

\vspace{3pt}

\noindent $(3)$ The singular locus of $\overline{\Sigma}^{*}$ coincides
with $\mP(0\oplus(U^{4})^{*})$, which is the $\psi$- image of $S_{\psi}$.

\vspace{3pt}

\noindent $(4)$ $\overline{\Sigma}^{*}$is a $8$-dimensional Fano
variety of degree $16$ with only Gorenstein canonical singularity
and with $-K_{\overline{\Sigma}^{*}}=\sO_{\overline{\Sigma}^{*}}(6)$.

\vspace{3pt}

\noindent $(5)$ The linear projection of $\mP((U^{3})^{*}\otimes U^{4}\oplus(U^{4})^{*})$
from $\mP(0\oplus(U^{4})^{*})$ induces the rational map $\overline{\Sigma}^{*}\dashrightarrow\mP((U^{3})^{*})\times\mP(U^{4})\simeq\mP^{2}\times\mP^{3}$.
\end{prop}

\begin{proof}
We can show the assertions in a quite similar way to the proof of
Proposition \ref{prop:DualG4}, so we only show (1). Let $[U^{1}]\in\mP(U^{4})$
be a point. The $\sigma$-fiber over the point $[U^{1}]$ is $\mP\left((U^{3})^{*}\otimes U^{1}\oplus(U^{4}/U^{1})^{*}\right).$
With this description, we immediately see that $\overline{\Sigma}^{*}$
is contained in the variety defined by the equation (\ref{eq:G9DualEq}),
which we temporarily denote by $(\overline{\Sigma}^{*})'$. Let $[D,\bm{p}]$
be a point of $(\overline{\Sigma}^{*})'$ such that $D\not=0$. Then,
by Lemma \ref{lem:ABFib} (2), we see that the $\psi$-fiber over
$[D,\bm{p}]$ consists of one point $[U^{1}]\in\mP(U^{4})$ such that
$U^{1}$ spans the image of the linear map $U^{3}\to U^{4}$ defined
by the rank 1 matrix $D$. Therefore, $\mP(\sE^{\perp})\to(\overline{\Sigma}^{*})'$
is dominant, hence is surjective, and is also birational. 
\end{proof}

\subsection{Curve of genus $9$ \label{subsec:Curve-of-genus9}}

In this subsection, we characterize smooth 1-dimensional linear sections
of $\overline{\Sigma}^{*}$ in the framework of the classification
of algebraic curves. We denote by $\pi_{1}$ the first projection
$\mP^{2}\times\mP^{3}\to\mP^{2}$ and by $\pi_{2}$ the second projection
$\mP^{2}\times\mP^{3}\to\mP^{3}$.

The following lemma should be well-known for experts but we include
a proof since we cannot find any reference.
\begin{lem}
\label{lem:notrisec} Let $\mP^{11}$ be the ambient space of the
Segre embedded $\mP^{2}\times\mP^{3}$. The following assertions hold:

\vspace{3pt}
\end{lem}

\noindent $(1)$ For a line $l$ in $\mP^{11}$, it holds that $l\subset\mP^{2}\times\mP^{3}$
or $l\cap(\mP^{2}\times\mP^{3})$ consists of at most two points.

\vspace{3pt}

\noindent $(2)$ Let $P$ be a plane in $\mP^{11}$. If $P\cap(\mP^{2}\times\mP^{3})$
contains infinite number of points, then $P\subset\mP^{2}\times\mP^{3}$,
or $P\cap(\mP^{2}\times\mP^{3})$ is a conic, a line or the union
of a line and a point. Otherwise, $P\cap(\mP^{2}\times\mP^{3})$ consists
of at most three points.
\begin{proof}
The assertion (1) immediately follows since $\mP^{2}\times\mP^{3}$
is defined by the quadrics. 

We show the assertion (2). The first assertion follows since $\mP^{2}\times\mP^{3}$
is defined by the quadrics. Therefore, for the second assertion, we
may assume that $P\cap(\mP^{2}\times\mP^{3})$ consists of a finite
number of points, and contains at least 3 points, say, $p_{1},p_{2},p_{3}$.
If $p_{1},p_{2},p_{3}$ are colinear, then $P\cap(\mP^{2}\times\mP^{3})$
contains the line they span by (1), a contradiction. Thus $p_{1},p_{2},p_{3}$
span the plane $P$. If two of them are contained in a fiber of $\pi_{1}$
or $\pi_{2}$, then $P\cap(\mP^{2}\times\mP^{3})$ contains the line
joining the two points, a contradiction. Thus no two of them are contained
in a fiber of $\pi_{1}$ or $\pi_{2}$. Let $q_{i}$ and $r_{i}$
be the images of $p_{i}$ by $\pi_{1}$ and $\pi_{2}$ respectively
($i=1,2,3$). Then $L_{q}:=\langle q_{1},q_{2},q_{3}\rangle$ and
$L_{r}:=\langle r_{1},r_{2},r_{3}\rangle$ are lines or planes. Note
that $P$ is contained in the ambient space of $L_{q}\times L_{r}$
and hence $P\cap(\mP^{2}\times\mP^{3})=P\cap(L_{q}\times L_{r})$.
If $L_{q}$ and $L_{r}$ are lines, then $P\cap(L_{q}\times L_{r})$
is a conic, a contradiction. If one of $L_{q}$ and $L_{r}$ is a
line and another is a plane, then $P\cap(L_{q}\times L_{r})$ consists
of three points $p_{1},p_{2},p_{3}$ as desired since $\deg(L_{q}\times L_{r})=3$.
Finally, assume that $L_{q}$ and $L_{r}$ are planes. Then, by coordinate
changes of $L_{q}$ and $L_{r}$ if necessary, we may assume that
$p_{1}=q_{1}=(1:0:0)$, $p_{2}=q_{2}=(0:1:0)$, and $p_{3}=q_{3}=(0:0:1)$.
Then it is easy to check that $P\cap(L_{q}\times L_{r})$ consists
of three points $p_{1},p_{2},p_{3}$ as desired.
\end{proof}
\begin{cor}
\label{cor:nontri} Let $C$ be a smooth non-hyperelliptic curve of
genus $9$. Assume that 

\vspace{3pt}

\noindent $(1)$ the canonical model of $C$ is contained in the
Segre embedded $\mP^{2}\times\mP^{3}$ (then we identify $C$ with
its canonical model), and 

\vspace{3pt}

\noindent $(2)$ the first projection $\pi_{1}$ induces a birational
map from $C$ onto a septic plane curve with only double points as
its singularities.

\vspace{3pt}

\noindent $(3)$ the second projection $\pi_{2}$ induces a birational
map from $C$ onto the curve with at worst double points as its singularities.

\vspace{3pt}

\noindent Then $C$ is not 4-gonal and the Clifford index of $C$
is $3$.
\end{cor}

\begin{proof}
Assume by contradiction that $C$ has a $g_{3}^{1}$, say, $\delta$.
Since $C$ is non-hyperelliptic, $|\delta$| has no base points. Let
$D\in|\delta|$ be a general element. By the Riemann-Roch theorem,
we have $h^{0}(K_{C}-D)=7,$ which implies that ${\rm Supp}\,D$ spans
a line $l_{D}$. By Lemma \ref{lem:notrisec} (1) and the assumption
(1), $l_{D}$ must be contained in $\mP^{2}\times\mP^{3}$. It is
easy to see that $l_{D}$ is contained in a fiber of $\pi_{1}$ or
$\pi_{2}$. Then the image by $\pi_{1}$ or $\pi_{2}$ of $C$ has
a triple point, a contradiction to the assumption (2) or (3). 

Assume by contradiction that $C$ has a $g_{4}^{1}$, say, $\varepsilon$.
Since $C$ is non-trigonal as we have seen above, $|\varepsilon$|
has no base points. Let $E\in|\varepsilon|$ be a general element.
By the Riemann-Roch theorem, we have $h^{0}(K_{C}-E)=6,$ which implies
that ${\rm Supp}\,E$ spans a plane $P_{E}$. By Lemma \ref{lem:notrisec}
(2) and the assumption (1), it holds that $P_{E}\cap(\mP^{2}\times\mP^{3})$
is $P_{E},$ or contains a line, or coincides with a smooth conic,
say, $q_{E}$. In the the first case, $P_{E}$ must be contained in
a fiber of $\pi_{1}$ or $\pi_{2}$, and then the image by $\pi_{1}$
or $\pi_{2}$ of $C$ has a quadruple point, a contradiction to the
assumption (2) or (3). In the the second case, at least three points
of $E$ is contained in the line. Since the line is contained in a
fiber of $\pi_{1}$ or $\pi_{2}$, the image by $\pi_{1}$ or $\pi_{2}$
of $C$ has a triple or a quadruple point, a contradiction to the
assumption (2) or (3). Assume the third case occurs. If the smooth
conic $q_{E}$ is contained in a fiber of $\pi_{1}$ or $\pi_{2}$,
we may derive a contradiction in the same way as the first and the
second cases. Therefore $q_{E}$ is mapped to a line by $\pi_{1}$
and $\pi_{2}$. Let $l_{E}:=\pi_{1}(q_{E})$. Then $\pi_{1}^{*}(l_{E_{1}})\cap C$
and $\pi_{1}^{*}(l_{E_{2}})\cap C$ are linearly equivalent for $E_{1},E_{2}\in|\varepsilon|$.
Since $E_{1}\subset\pi_{1}^{*}(l_{E_{1}})\cap C$ and $E_{2}\subset\pi_{1}^{*}(l_{E_{2}})\cap C$
and $E_{1}\sim E_{2}$, it holds that $(\pi_{1}^{*}(l_{E_{1}})\cap C)-E_{1}$
and $(\pi_{1}^{*}(l_{E_{2}})\cap C)-E_{2}$ are linearly equivalent.
Since $\deg(\pi_{1}^{*}(l_{E_{i}})\cap C-E_{i})=3$, this implies
that $C$ is trigonal, a contradiction. Therefore we have shown that
$C$ is not 4-gonal. 

Now the assertion that the Clifford index of $C$ is $3$ follows
from \cite[Cor. 3.1.1 and 3.2.1]{Sa}.
\end{proof}
The following result for a curve of genus 9 is similar to the one
for a curve of genus 8 as in \cite{MuId}.
\begin{prop}
\label{prop:curveG9V1} Let $C$ be a smooth curve of genus $9$.
The following assertions $(1)$ and $(2)$ are equivalent:

\vspace{3pt}

\noindent $(1)$ There exists a birational morphism $\iota_{1}$
from $C$ to a septic plane curve $C_{1}$ with only double points
and an isomorphism $\iota_{2}\colon C\to C_{2}$ to a space curve
$C_{2}$ of degree $9$ such that $\iota_{1}^{*}\sO_{C_{1}}(1)+\iota_{2}^{*}\sO_{C_{2}}(1)=K_{C}$. 

\vspace{3pt}

\noindent $(2)$ $C$ is isomorphic to the complete intersection
in $\mP^{2}\times\mP^{3}$ of three divisors of $(1,1)$-type and
a divisor of $(1,2)$-type. 
\end{prop}

\begin{proof}
\vspace{3pt}

\noindent $(1)\Rightarrow(2)$. Assume that the assertion (1) holds.
It is easy to check by standard calculations for curves on surfaces
that the curve $C_{2}$ is not contained in a plane nor a quadric
surface since $\deg C_{2}=9$ and $g(C_{2})=9.$ Let $\iota_{3}\colon C\to\mP^{2}\times\mP^{3}\hookrightarrow\mP^{11}$
be the composite of the morphism induced by $\iota_{1}\times\iota_{2}$,
and the Segre embedding. Since $\iota_{2}$ is an isomorphism, $\iota_{3}$
defines an isomorphism onto the image, which we denote by $C_{3}$.
Note that, by the construction and the condition that $\iota_{1}^{*}\sO_{C_{1}}(1)+\iota_{2}^{*}\sO_{C_{2}}(1)=K_{C}$,
the restriction to $C_{3}$ of the divisor of $(1,1)$-type in $\mP^{2}\times\mP^{3}$
is $K_{C_{3}}.$ By the Riemann-Roch theorem, we see that there are
at least three linearly independent forms of bidegree $(1,1)$ on
$\mP^{2}\times\mP^{3}$ vanishing on $C_{3}$. We take any such three
$\eta_{1},\eta_{2},\eta_{3}$. Let $x_{1},x_{2},x_{3}$ be coordinates
of $\mP^{2}$ and $y_{1},\dots,y_{4}$ coordinates of $\mP^{3}$.
We may write 
\[
\left(\begin{array}{ccc}
\eta_{1} & \eta_{2} & \eta_{3}\end{array}\right)=\left(\begin{array}{ccc}
x_{1} & x_{2} & x_{3}\end{array}\right)M,
\]
where $M$ is a certain $3\times3$ matrix whose entries are linear
forms with respect to $y_{1},\dots,y_{4}$. 

We show that $\det M$ is a nonzero cubic form. Assume by contradiction
that $\det M\equiv0$. Then we may consider $M$ defines a 3-dimensional
linear subspace in the determinantal cubic hypersurface of the generic
$3\times3$ matrix. By the classification result in \cite{At} (see
also \cite[Thm.1.1]{EH}, \cite[7A]{CI}), we have the following possibilities
of $M$ by changing coordinates of $\mP^{2}$ and $\mP^{3}$, and
$\eta_{1},\eta_{2},\eta_{3}$ if necessary:
\[
\left(\begin{array}{ccc}
0 & * & *\\
0 & * & *\\
0 & * & *
\end{array}\right),\left(\begin{array}{ccc}
0 & 0 & *\\
0 & 0 & *\\
* & * & *
\end{array}\right),\left(\begin{array}{ccc}
* & * & *\\
0 & 0 & *\\
0 & 0 & *
\end{array}\right),\left(\begin{array}{ccc}
* & * & *\\
* & * & *\\
0 & 0 & 0
\end{array}\right).
\]
The first case is impossible since then $\eta_{1}=0.$ In the second
and third cases, each of $\eta_{1},\eta_{2}$ is the product of linear
forms of $\mP^{2}$ and $\mP^{3}$. Then $C_{1}$ is a line or $C_{2}$
is a plane curve, a contradiction. Now we consider the 4th case. We
write the matrix more explicitly as $M=\left(\begin{array}{ccc}
a_{1} & a_{2} & a_{3}\\
b_{1} & b_{2} & b_{3}\\
0 & 0 & 0
\end{array}\right),$ where $a_{i},b_{j}$ are linear forms of $y_{1},\dots,y_{4}$. Then
the locus $\{\eta_{1}=\eta_{2}=\eta_{3}=0\}$ contains $\{(0:0:1)\}\times\mP^{3}$
and the image in $\mP^{3}$ of $\{\eta_{1}=\eta_{2}=\eta_{3}=0\}\setminus\{(0:0:1)\}\times\mP^{3}$
is contained in the locus $S:=\left\{ \rank\left(\begin{array}{ccc}
a_{1} & a_{2} & a_{3}\\
b_{1} & b_{2} & b_{3}
\end{array}\right)\leq1\right\} $. Since $C_{3}$ cannot be contained in $\{(0:0:1)\}\times\mP^{3}$,
we see that $C_{2}\subset S$. Since $S$ cannot coincide with $\mP^{3}$
and defined by quadrics, $C_{2}$ is contained in a quadric surface,
a contradiction. Therefore we have shown that $\det M$ is a nonzero
cubic form for any choice of $\eta_{1},\eta_{2},\eta_{3}$. Moreover,
we see that $C_{2}\subset\{\det M=0\}$. Since $C_{2}$ is not contained
in a plane nor a quadric, the cubic surface $\{\det M=0\}$ is irreducible.

Now we show that there are \textit{exactly} three linearly independent
forms of bidegree $(1,1)$ on $\mP^{2}\times\mP^{3}$ vanishing on
$C_{3}$. Assume the contrary. Let $\eta_{1},\eta_{2},\eta_{3},\eta_{4}$
be four linearly independent forms of bidegree $(1,1)$ on $\mP^{2}\times\mP^{3}$
vanishing on $C_{3}$. Let $M_{1}$ and $M_{2}$ be the matrix $M$
defined for $\eta_{1},\eta_{2},\eta_{3}$ and $\eta_{1},\eta_{2},\eta_{4}$
respectively. We have shown that $C_{2}\subset\{\det M_{1}=0\}\cap\{\det M_{2}=0\}.$
If $\{\det M_{1}=0\}\not=\{\det M_{2}=0\}$, then $C_{2}=\{\det M_{1}=0\}\cap\{\det M_{2}=0\}$
by the reason of degree since $\{\det M_{1}=0\}$ and $\{\det M_{2}=0\}$
are irreducible. However, the genus of $C_{2}$ is not equal to that
of the curve $\{\det M_{1}=0\}\cap\{\det M_{2}=0\}$, a contradiction.
Therefore we must have \{$\det M_{1}=0\}=\{\det M_{2}=0\}$. Then,
for suitable nonzero constants $\alpha,\beta$, we have $\det(\alpha M_{1}+\beta M_{2})\equiv0$.
By the previous paragraph, this implies that $\eta_{1},\eta_{2},\alpha\eta_{3}+\beta\eta_{4}$
must be linearly dependent, a contradiction. Thus we have shown that
there are exactly three linearly independent forms of bidegree $(1,1)$
on $\mP^{2}\times\mP^{3}$ vanishing on $C_{3}$, and hence $C_{3}$
is a canonical curve. 

In particular, $C_{3}$ is non-hyperelliptic and then $H^{0}(K_{C_{3}})$
generates the canonical ring of $C_{3}$ by the Max Noether theorem
(cf.~\cite[p.117]{ACGH}). Therefore, by the Riemann-Roch theorem,
we see that the number of linearly independent quadratic forms on
$\langle C_{3}\rangle$ vanishing along $C_{3}$ is $45-24=21$. 

Since $\langle C_{3}\rangle$ is of codimension 3 in $\mP^{11}$,
the dimension of the space of forms of bidegree $(1,2)$ in $\mP^{2}\times\mP^{3}\cap\langle C_{3}\rangle$
is at least $3\times10-3\times4=18$. On the other hand, since $\deg\sO_{C_{3}}(1,2)=25$,
we have $h^{0}(\sO_{C_{3}}(1,2))=17$ by the Riemann-Roch theorem.
Therefore there exists a form $\xi$ of bidegree $(1,2)$ in $\mP^{2}\times\mP^{3}\cap\langle C_{3}\rangle$
vanishing along $C_{3}$. Producting $\xi$ with three linearly independent
forms $x_{1},x_{2},x_{3}$ of bidegree $(1,0)$, we obtain three linearly
independent quadratic forms on $\langle C_{3}\rangle$ vanishing along
$C_{3}$. These and the 18 quadratic forms defining $\mP^{2}\times\mP^{3}\cap\langle C_{3}\rangle$
are clearly linearly independent, thus they form the 21-dimensional
space of quadratic forms on $\langle C_{3}\rangle$ vanishing along
$C_{3}$. Since $C_{3}$ satisfies the assumptions of Corollary \ref{cor:nontri},
$C_{3}$ is non-trigonal. Therefore $C_{3}$ is scheme theoretically
the intersection of quadrics containing $C_{3}$ by Enriques-Babbage-Petri
theorem (cf. \cite[p.124,p.131]{ACGH}). Thus it holds that $C_{3}=\left(\mP^{2}\times\mP^{3}\cap\langle C_{3}\rangle\right)\cap\left\{ \xi=0\right\} $
scheme-theoretically, and hence the assertion (2) holds.

\vspace{3pt}

\noindent $(2)\Rightarrow(1)$. Assume that the assertion (2) holds.
We identify $C$ with its model in $\mP^{2}\times\mP^{3}$ as in (2).
Let $S_{C}$ be the complete intersection of three divisors of type
$(1,1)$ containing $C$. Since $C$ is a smooth curve and an ample
divisor on $S_{C}$, we see that $S_{C}$ is an irreducible Gorenstein
surface and is smooth along $C$. Since a divisor of type $(1,1)$
is degree 1 on a fiber of $\pi_{1}$ or $\pi_{2}$, the restriction
of a fiber of $\pi_{1}$ or $\pi_{2}$ to $S_{C}$ is a point or a
line. Then we see that $S_{C}$ is rational since $S_{C}$ birationally
dominates $\mP^{2}$, and $\pi_{2}(S_{C})$ is a cubic surface since
it follows that $\deg\sO_{S_{C}}(0,1)=3$. 

We show that $S_{C}$ is normal. Since $S_{C}\to\mP^{2}$ is birational
and $\mP^{2}$ is smooth, $S_{C}$ is possibly non-normal only along
lines which are the restrictions of $\pi_{1}$-fibers to $S_{C}$
by the Zariski main theorem. However, $S_{C}$ is smooth along $C$
and such lines intersects $C$ since $C$ is an ample divisor on $S_{C}$,
$S_{C}$ cannot be non-normal.

We check the desired properties of $\pi_{1}(C)$. By the assumption
(2), it is easy to obtain $\deg\sO_{C}(1,0)=7$. Let $l$ be a line
(if it exists) which is the restriction of a $\pi_{1}$-fiber to $S_{C}$.
Since $C$ is the restriction to $S_{C}$ of a divisor of type $(1,2)$,
we have $C\cdot l=2$. Since it holds that $-K_{S_{C}}\cdot l=1$,
we see that $\pi_{1}(C)$ has a double point at $\pi_{1}(l)$ by \cite[Thm.0.1]{LS}.
Therefore $\pi_{1}$ induces a birational map from $C$ to a septic
plane curve with only double points as singularities. 

We check the desired properties of $\pi_{2}(C)$. By the assumption
(2), it is easy to obtain $\deg\sO_{C}(0,1)=9$. Let $m$ be a line
which is the restriction of a $\pi_{2}$-fiber to $S_{C}$. Since
$C$ is the restriction to $S_{C}$ of a divisor of type $(1,2)$,
we have $C\cdot m=1$. Since it holds that $-K_{S_{C}}\cdot m=1$,
we see that $\pi_{2}(S_{C})$ and $\pi_{2}(C)$ are smooth at $\pi_{1}(m)$
by \cite[Thm.0.1]{LS}. Now let $p$ be a point of $C$ such that
$S_{C}\to\pi_{2}(S_{C})$ is finite near $p$. We may choose two divisors
$D_{1},D_{2}$ of $(1,1)$-type containing $S_{C}$ such that $D_{1}\cap D_{2}\to\mP^{3}$
is finite near $p$. Then, $D_{1}\cap D_{2}$ is a section of the
$\mP^{2}$-bundle $\mP^{2}\times\mP^{3}\to\mP^{3}$ near $p$, hence
$D_{1}\cap D_{2}\to\mP^{3}$ is an isomorphism near $p$. Therefore
$S_{C}\to\pi_{2}(S_{C})$ is also an isomorphism near $p$. Now we
have seen that $C\to\pi_{2}(C)$ is isomorphic at any point as desired. 

Thus we have verified that the assertion (2) holds.
\end{proof}
\begin{rem}
\label{rem:G9G27}\vspace{3pt}

\noindent (1) In the proof of $(1)\Rightarrow(2)$ of Proposition
\ref{prop:curveG9V1}, a famous classic construction of a cubic surface
by C.$\,$Segre \cite{Se} (see also \cite[Sect.2]{D}) naturally
appears.
\vspace{3pt}

\noindent (2) The assumption on $C_{1}$ as in Proposition \ref{prop:curveG9V1}
is natural in view of gonality and Clifford index (cf.~\cite{Sa}).
The assumption on $C_{2}$, however, is more delicate as we see in
the following example. 

Let $C_{1}\subset\mP^{2}$ be a septic 6-nodal plane curve such that
the 6 nodes are located on a smooth conic $q$. Let $S_{C}\to\mP^{2}$
be the blow-up at the six nodes of $C_{1}$ and $f_{i}\,(1\leq i\leq6)$
the exceptional curves. By the assumption, $S_{C}$ is a cubic weak
del Pezzo surface and there is a birational morphism from $S_{C}$
to a cubic surface $T$ contracting the strict transform $q'$ of
$q$. The strict transform $\widetilde{C}_{2}$ of $C_{1}$ on $S_{C}$
is smooth and is linearly equivalent to $7m-2\sum_{i=1}^{6}f_{i}$
and $q'\sim2m-\sum_{i=1}^{6}f_{i}$, where $m$ is the total transform
on $S_{C}$ of a line of $\mP^{2}$. In this case, we have a naturally
induced morphism $\widetilde{C}_{2}\to\mP^{2}\times\mP^{3}$ and it
induces an isomorphism from $\widetilde{C}_{2}$ onto the image $C_{3}$
since $\widetilde{C}_{2}\to C_{3}$ is an isomorphism outside $q'\cap\widetilde{C}_{2}$
and $\widetilde{C_{2}}\to C_{1}$ is isomorphism near $q'\cap\widetilde{C}_{2}$.
Since $\widetilde{C}_{2}\cdot q'=2$, the image $C_{2}\subset T$
of $\widetilde{C}_{2}$ has a double point at the image of $q'$ as
its singularity. Therefore, by Proposition \ref{prop:curveG9V1},
$C_{3}$ cannot be the complete intersection in $\mP^{2}\times\mP^{3}$
of three divisors of $(1,1)$-type and a divisor of $(1,2)$-type.
However, the Clifford index of $C_{3}$ is 3 by Corollary \ref{cor:nontri}.
Note that, since $C_{3}\sim3(3m-\sum_{i=1}^{6}f_{i})-q'$, and $3m-\sum_{i=1}^{6}f_{i}$
is the pull-back of $\sO_{T}(1)$, we see that $C_{2}$ is the complete
intersection between $T$ and another cubic surface $T'$. 
\end{rem}

The following result connects a property of curve of genus 9 with
the key variety.
\begin{cor}
\label{cor:curveG9V2} Let $C$ be a smooth curve of genus $9$. The
following assertions $({\rm a})$ and $({\rm b})$ are equivalent: 

\vspace{3pt}
\end{cor}

\noindent $({\rm a})$ There exists a birational morphism $\iota_{1}$
from $C$ to a septic plane curve $C_{1}$ with only double points
and an isomorphism $\iota_{2}\colon C\to C_{2}$ to a space curve
$C_{2}$ of degree $9$ such that $\iota_{1}^{*}\sO_{C_{1}}(1)+\iota_{2}^{*}\sO_{C_{2}}(1)=K_{C}.$ 

\vspace{3pt}

\noindent $({\rm b})$ $C$ is isomorphic to a linear section of
$\overline{\Sigma}^{*}$.
\begin{proof}
It suffices to show the equivalence of (b) and the assertion (2) of
Proposition \ref{prop:curveG9V1}.

Assume that the assertion (2) of Proposition \ref{prop:curveG9V1}
holds. Let $C_{3}$ and $\xi$ be as in the proof of Proposition \ref{prop:curveG9V1}
$(1)\Rightarrow(2)$, and $\widetilde{\xi}$ a form of bidegree $(1,2)$
on $\mP^{2}\times\mP^{3}$ which is a lift of $\xi$. We write $\widetilde{\xi}=l_{1}y_{1}+\cdots+l_{4}y_{4}$,
where $l_{i}$ are forms of bidegree $(1,1)$. Let $L_{i}$ be the
linear forms on $\mP^{11}$ corresponding to $l_{i}$. We have 
\[
\{\widetilde{\xi}=0\}=\mP^{2}\times\mP^{3}\cap\left\{ \left(\begin{array}{cccc}
r_{11} & r_{12} & r_{13} & r_{14}\\
r_{21} & r_{22} & r_{23} & r_{24}\\
r_{31} & r_{32} & r_{33} & r_{34}
\end{array}\right)\left(\begin{array}{c}
L_{1}\\
L_{2}\\
L_{3}\\
L_{4}
\end{array}\right)=\bm{o}\right\} ,
\]
where $r_{ij}$ are coordinates of $\mP^{11}$. We consider $r_{ij}$
are the entries of $\empty^{t}\!D$ and let $p_{1},\dots,p_{4}$ be
the entries of $\bm{p}$ as in the equations (\ref{eq:G9DualEq})
of $\overline{\Sigma}^{*}$. Then we see that $\{\widetilde{\xi}=0\}$
is projectively equivalent to $\overline{\Sigma}^{*}\cap\left\{ p_{1}=L_{1},p_{2}=L_{2},p_{3}=L_{3},p_{4}=L_{4}\right\} $,
which is a linear section of $\overline{\Sigma}^{*}$. Therefore,
finally, we have seen that $C_{3}$ is projectively equivalent to
a linear section of $\overline{\Sigma}^{*}$ as desired.

The converse follows by reversing the above discussion.
\end{proof}
\begin{rem}
We do not see the relation of $\overline{\Sigma}^{*}$ with the symplectic
Grassmanian ${\rm Sp}(3,6)$.
\end{rem}

\section{\textbf{$\mQ$-Fano 3-fold of Genus 6 and $\mathsf{\mathtt{Q}}$-type
\label{sec:genus 6Q}}}

In this case, it holds that $\mP(\sE)\simeq\mP(\sE^{\perp})$ since
$\sU|_{Q^{3}}\simeq\sQ^{*}|_{Q^{3}}$. This self-duality could be
compared to that of the orthogonal Grassmanian ${\rm OG(5,10)}$ (see
\cite{Mu4}). We only give a few remark about a curve $C$ of genus
6 which is a linear section of $\overline{\Sigma}^{*}$. By self-duality,
we identify $\overline{\Sigma}^{*}$with $\overline{\Sigma}$. Since
$C$ is a smooth linear section of $\overline{\Sigma}$, we see that
$C$ is disjoint from the singular locus of $\overline{\Sigma}.$
Therefore $C$ can be consider as a linear section of a quadric section
of $A_{\mathtt{Q}}$. It is easy to see the converse holds; if $C$
is a linear section of a quadric section of $A_{\mathtt{Q}}$, then
$C$ is also a linear section of $\overline{\Sigma}$. By \cite[Sect. 5]{Mu3},
a general smooth curve of genus 6 is a linear section of a quadric
section of $A_{\mathtt{Q}}$, hence is a linear section of $\overline{\Sigma}$.


\section{\textbf{$\mQ$-Fano 3-fold of Genus 6 and $\mathsf{\mathtt{C}}$-type}}

\subsection{Descriptions of $\mP(\sE^{\perp})$ \label{subsec:DescriptionsPEperpG6C}}

We use the notation in Table 1. Let $y_{1},\dots,y_{5}$ be coordinates
of $U^{5}$ and we may assume that the twisted cubic $\gamma_{\mathtt{C}}$
is equal to 
\[
\left\{ y_{2}^{2}-y_{1}y_{3}=y_{1}y_{4}-y_{2}y_{3}=y_{3}^{2}-y_{2}y_{4}=y_{5}=0\right\} .
\]

Let $F_{a}$ and $F_{b}$ be the $a$- and $b$-exceptional divisors,
respectively. Since $b\circ a^{-1}\colon A_{\mathtt{C}}\dashrightarrow\mP(U^{5})$
is the restriction of the linear projection from $\Pi$, we have 
\begin{equation}
b^{*}\sO_{\mP(U^{5})}(1)=a^{*}\sO_{A_{\mathtt{C}}}(1)-F_{a}.\label{eq:LBLA}
\end{equation}

As in \cite[the subsec.2.4]{Tak2}, we consider $A_{\mathtt{C}}\subset\mP(\wedge^{2}V^{3}\oplus U^{5})$,
where $V^{3}$ is a 3-dimensional vector space.
\begin{prop}
\label{prop:Dual Cubic} The following assertions hold:

\vspace{3pt}

\noindent $(1)$ $\overline{\Sigma}^{*}$ is the cubic hypersurface
in $\mP((\wedge^{2}V^{3})^{*}\oplus(U^{5})^{*}\oplus U^{5})\simeq\mP^{12}$
with the following equation: 
\begin{equation}
p_{1}(q_{2}^{2}-q_{1}q_{3})+p_{2}(q_{1}q_{4}-q_{2}q_{3})+p_{3}(q_{3}^{2}-q_{2}q_{4})+q_{5}\left(\sum_{i=1}^{5}r_{i}q_{i}\right)=0,\label{eq:CubicEq}
\end{equation}
where $q_{1},\dots,q_{5}$ are coordinates of $U^{5}$, $p_{1},p_{2},p_{3}$
are those of $(\wedge^{2}V^{3})^{*}$, and $r_{1},\dots,r_{5}$ are
those of $(U^{5})^{*}$. The morphism $\psi$ is birational onto $\overline{\Sigma}^{*}$.

\vspace{3pt}

\noindent $(2)$ The singular locus of $\overline{\Sigma}^{*}$ is
the union of $\mP((\wedge^{2}V^{3})^{*}\oplus(U^{5})^{*}\oplus0)\simeq\mP^{7}$
and the closure $S_{F}$ of the $6$-dimensional locus 
\[
\left\{ q_{1}=1,q_{3}=q_{2}^{2},q_{4}=q_{2}^{3},q_{5}=0,p_{1}=q_{2}^{2}p_{3},p_{2}=q_{2}p_{3},r_{1}=-q_{2}r_{2}-q_{2}^{2}r_{3}-q_{2}^{3}r_{4}\right\} .
\]
The cubic $\overline{\Sigma}^{*}$ has ordinary double points generically
along each of the irreducible components of $\Sing\overline{\Sigma}^{*}$.

\vspace{3pt}

\noindent $(3)$ The $\psi$-exceptional locus is the union of the
two divisors
\[
E_{\mP(\sE^{\perp})}:=\mP(a^{*}(\Omega_{\mP^{7}}^{1}(1))\oplus0),\text{{and}}\,\,F_{\mP(\sE^{\perp})}:=\sigma^{*}F_{b},
\]
and $\psi(E_{\mP(\sE^{\perp})})=\mP((\wedge^{2}V^{3})^{*}\oplus(U^{5})^{*}\oplus0)\simeq\mP^{7}$
and $\psi(F_{\mP(\sE^{\perp})})=$ 
\[
\{q_{2}^{2}-q_{1}q_{3}=q_{1}q_{4}-q_{2}q_{3}=q_{3}^{2}-q_{2}q_{4}=q_{5}=0\},
\]
 where the latter is the cone over a twisted cubic with $\mP((\wedge^{2}V^{3})^{*}\oplus(U^{5})^{*}\oplus0)$
as the vertex. The image of the $\psi$-exceptional locus contains
$\Sing\overline{\Sigma}^{*}$. 

\vspace{3pt}

\noindent $(4)$ The morphism $\psi$ is the blow-up along $\psi(F_{\mP(\sE^{\perp})})$
outside $\Sing\overline{\Sigma}^{*}.$ The $\psi$-fiber over a point
$\mathsf{t}\in\mP((\wedge^{2}V^{3})^{*}\oplus(U^{5})^{*}\oplus0)$
is isomorphic to the total transform of the hyperplane section of
$A_{\mathtt{C}}$ corresponding to $\mathsf{t}$ by projective duality.
In particular, the $\psi$-fiber over a general point $\mathsf{\mathsf{t}}\in\mP((\wedge^{2}V^{3})^{*}\oplus(U^{5})^{*}\oplus0)$
is the smooth 3-fold obtained by blowing up $B_{5}$ along a line.
\end{prop}

\begin{proof}
(1). First we see that the morphism $\psi$ is birational onto a certain
cubic hypersurface computing $H_{\mP(\sE^{\perp})}^{11}$. It is well-known
that
\begin{align*}
c_{t}(T_{\mP(\wedge^{2}V^{3}\oplus U^{5})}(-1)|_{A_{\mathtt{C}}}) & =1+c_{1}(\sO_{A_{\mathtt{C}}}(1))t+c_{1}(\sO_{A_{\mathtt{C}}}(1))^{2}t^{2}+c_{1}(\sO_{A_{\mathtt{C}}}(1))^{3}t^{3}+c_{1}(\sO_{A_{\mathtt{C}}}(1))^{4}t^{4}\\
 & =1+c_{1}(\sO_{A_{\mathtt{C}}}(1))t+c_{1}(\sO_{A_{\mathtt{C}}}(1))^{2}t^{2}+(5l)t^{3}+5t^{4},
\end{align*}
where $l$ is the class of a line in $A_{\mathtt{C}}.$ We set $c_{A}:=c_{1}(a^{*}\sO_{A_{\mathtt{C}}}(1))$
and $c_{B}:=c_{1}(b^{*}\sO_{\mP(U^{5})}(1))$. By a standard computation,
we have 
\begin{align*}
c_{t}((\sE^{\perp})^{*}) & =1+(c_{A}+c_{B})t+(c_{A}c_{B}+c_{A}^{2})t^{2}+(c_{A}^{2}c_{B}+a^{*}(5l))t^{3}+\left(c_{B}.\left(a^{*}(5l)\right)+5\right)t^{4}.
\end{align*}
From this, we have 
\begin{align*}
H_{\mP(\sE^{\perp})}^{11}= & s_{4}((\sE^{\perp})^{*})\\
= & c_{1}((\sE^{\perp})^{*})^{4}-3c_{1}((\sE^{\perp})^{*})^{2}c_{2}((\sE^{\perp})^{*})+2c_{1}((\sE^{\perp})^{*})c_{3}((\sE^{\perp})^{*})+c_{2}((\sE^{\perp})^{*})^{2}-c_{4}((\sE^{\perp})^{*})\\
= & 6-c_{A}^{3}c_{B}+c_{A}c_{B}^{3}.
\end{align*}
By (\ref{eq:LBLA}), we have 
\begin{equation}
c_{B}=c_{A}-c_{1}(F_{a}).\label{eq:cBcA}
\end{equation}
By \cite[Sect.10]{Fuj}, $b$ is the blow-up of $\mP(U^{5})$ along
a twisted cubic curve and the $b$-exceptional divisor $F_{b}$ is
linearly equivalent to $a^{*}\sO_{A_{\mathtt{C}}}(1)-2F_{a}$. Therefore,
together with (\ref{eq:cBcA}), we have 
\begin{equation}
c_{A}=2c_{B}-c_{1}(F_{b}).\label{eq:cAcB}
\end{equation}
 By (\ref{eq:cBcA}) and (\ref{eq:cAcB}), we have $c_{A}^{3}c_{B}=c_{A}^{3}(c_{A}-c_{1}(F_{a}))=5$
and $c_{A}c_{B}^{3}=(2c_{B}-c_{1}(F_{b}))c_{B}^{3}=2$. Therefore,
we have $H_{\mP(\sE^{\perp})}^{11}=3$. This implies that the $\psi$
is a birational morphism onto a cubic hypersurface in $\mP^{12}$
or is generically a triple cover of $\mP^{11}.$ The latter, however,
is impossible since $h^{0}(H_{\widehat{\Sigma}{}^{*}})=13$. 

Now we show that the equation of the $\psi$-image can be taken as
(\ref{eq:CubicEq}). We choose the equation of the twisted cubic $\gamma_{\mathtt{C}}$
which is the center of $b\colon\widehat{A}_{\mathtt{C}}\to\mP(U^{5})$
as above. Note that, since $b\circ a^{-1}\colon A_{\mathtt{C}}\dashrightarrow\mP(U^{5})$
is the projection from $\Pi$, $\widehat{A}_{\mathtt{C}}$ is contained
in $A_{\mathtt{C}}\times\mP(U^{5})$. Let $\mathsf{v}:=[\bm{x}]\times[\bm{y}]\in\widehat{A}_{\mathtt{C}}$
be a point with $\bm{x}\in\wedge^{2}V^{3}\oplus U^{5}$ and $\bm{y}\in U^{5}$.
The fiber of $\mP(\sE^{\perp})$ at $\mathsf{v}$ is equal to $\mP\left(\left((\wedge^{2}V^{3}\oplus U^{5})/\mC\bm{x}\right)^{*}\oplus\mC\bm{y}\right)$.
Therefore, by Lemma \ref{lem:ABFib} (2), for a point $[\bm{p}+\bm{r}+\bm{q}]\in\mP((\wedge^{2}V^{3})^{*}\oplus(U^{5})^{*}\oplus U^{5})$
with $\bm{p}\in(\wedge^{2}V^{3})^{*},\bm{r}\in(U^{5})^{*}$ and $\bm{q}\in U^{5}$,
the $\psi$-fiber over $[\bm{p}+\bm{r}+\bm{q}]$ consists of $[\bm{x}]\times[\bm{y}]$
such that $\bm{p}+\bm{r}\in\left((\wedge^{2}V^{3}\oplus U^{5})/\mC\bm{x}\right)^{*}$
and $\bm{q}\in\mC\bm{y}.$ Now assume that $\bm{q}\not=0$ and $[\bm{q}]\not\in\gamma_{\mathtt{C}}$,
and such a point $[\bm{x}]\times[\bm{y}]\in A_{\mathtt{C}}$ exists.
Then, the condition that $\bm{q}\in\mC\bm{y}$ is equivalent to that
$\bm{\bm{y}}\in\mC\bm{q}$, and hence we may assume that $\bm{y}=\bm{q},$
which we write as $\empty^{t}\left(\begin{array}{ccc}
q_{1} & \cdots & q_{5}\end{array}\right).$ Then, by the equality $a^{*}\sO_{A_{\mathtt{C}}}(1)=b^{*}\sO_{\mP(U^{5})}(2)-F_{b}$
and the condition that $[\bm{q}]\not\in\gamma_{\mathtt{C}}$, we may
write 
\[
\bm{x}=\empty^{t}\left(\begin{array}{cccccccc}
q_{2}^{2}-q_{1}q_{3} & q_{1}q_{4}-q_{2}q_{3} & q_{3}^{2}-q_{2}q_{4} & q_{5}q_{1} & q_{5}q_{2} & q_{5}q_{3} & q_{5}q_{4} & q_{5}^{2}\end{array}\right)
\]
taking suitable coordinates of $A_{\mathtt{C}}.$ Taking the coordinates
$p_{1},p_{2},p_{3},r_{1},\dots,r_{5}$ of $(\wedge^{2}V^{3})^{*}\oplus(U^{5})^{*}$
dual to $\wedge^{2}V^{3}\oplus U^{5}$, we see that the condition
that $\bm{p}+\bm{r}\in\left((\wedge^{2}V^{3}\oplus U^{5})/\mC\bm{x}\right)^{*}$
is equivalent to the equation of the cubic as in (\ref{eq:CubicEq}).
Therefore, a point $[\bm{p}+\bm{r}+\bm{q}]$ in the $\psi$-image
with $\bm{q}\not=0$ and $[\bm{q}]\not\in\gamma_{\mathtt{C}}$ is
contained in the cubic (\ref{eq:CubicEq}). Since we have seen the
$\psi$-image is also a cubic, it must coincide with the cubic (\ref{eq:CubicEq}).

The assertion (2) follows from straightforward calculations, which
we omit. 

\vspace{3pt}

\noindent (3). Let $\delta$ be a $\psi$-exceptional curve. Note
that $H_{\mP(\sE^{\perp})}\cdot\delta=0$. Moreover, it holds that
either $\sigma^{*}a^{*}\sO_{A_{\mathtt{C}}}(1)\cdot\delta>0$ or $\sigma^{*}b^{*}\sO_{\mP(U^{5})}(1)\cdot\delta>0$
since $\rho(\mP(\sE^{\perp}))=3$ and both $\sigma^{*}a^{*}\sO_{A_{\mathtt{C}}}(1)$
and $\sigma^{*}b^{*}\sO_{\mP(U^{5})}(1)$ are nef. Since $E_{\mP(\sE^{\perp})}\sim H_{\mP(\sE^{\perp})}-\sigma^{*}b^{*}\sO_{\mP(U^{5})}(1)$,
we have $\delta\subset E_{\mP(\sE^{\perp})}$ if $\sigma^{*}b^{*}\sO_{\mP(U^{5})}(1)\cdot\delta>0$.
Since $F_{\mP(\sE^{\perp})}\sim\sigma^{*}(b^{*}\sO_{\mP(U^{5})}(2)-a^{*}\sO_{A_{\mathtt{C}}}(1))$
by (\ref{eq:LBLA}), we have $\delta\subset F_{\mP(\sE^{\perp})}$
if $\sigma^{*}b^{*}\sO_{\mP(U^{5})}(1)\cdot\delta=0$ and $\sigma^{*}a^{*}\sO_{A_{\mathtt{C}}}(1)\cdot\delta>0.$
Therefore the $\psi$-exceptional locus is contained in $E_{\mP(\sE^{\perp})}\cup F_{\mP(\sE^{\perp})}$.
By the construction, it is obvious that $\psi(E_{\mP(\sE^{\perp})})=\mP((\wedge^{2}V^{3})^{*}\oplus(U^{5})^{*}\oplus0).$
Therefore, $E_{\mP(\sE^{\perp})}$ is contained in the $\psi$-exceptional
locus since $\dim\psi(E_{\mP(\sE^{\perp})})<\dim E_{\mP(\sE^{\perp})}$.
Since $F_{b}$ is the exceptional divisor of the blow-up of $\mP(U^{5})$
along the twisted cubic $\gamma_{\mathtt{C}}$ with the equation as
above, we see that $\psi(F_{\mP(\sE^{\perp})})$ is defined by the
same equation in $\mP((\wedge^{2}V^{3})^{*}\oplus(U^{5})^{*}\oplus U^{5})$
by the descriptions of $\psi$-fibers as in the proof of (1). The
divisor $F_{\mP(\sE^{\perp})}$ is contained in the $\psi$-exceptional
locus since $\dim\psi(F_{\mP(\sE^{\perp})})<\dim F_{\mP(\sE^{\perp})}$.

By a straightforward calculation, we see that the image of the $\psi$-exceptional
locus contains $\Sing\overline{\Sigma}^{*}$. 

\vspace{3pt}

\noindent (4). We show the first assertion. Let $\delta$ be a $\psi$-exceptional
curve. By the proof of (3), $\psi(\delta)\in E_{\mP(\sE^{\perp})}\cup F_{\mP(\sE^{\perp})}$,
and if $\psi(\delta)\in E_{\mP(\sE^{\perp})}$, then $\psi(\delta)$
is a singular point of $\overline{\Sigma}^{*}$. From now on, we assume
that $\psi(\delta)\not\in E_{\mP(\sE^{\perp})}$ and $\psi(\delta)\in F_{\mP(\sE^{\perp})}$.
By the proof of (3) again, we have $\sigma^{*}b^{*}\sO_{\mP(U^{5})}(1)\cdot\delta=0$
and $\sigma^{*}a^{*}\sO_{A_{\mathtt{C}}}(1)\cdot\delta>0.$ Since
$H_{\mP(\sE^{\perp})}\cdot\delta=0$, $\delta$ cannot be contracted
by $\sigma$. Therefore, by $\sigma^{*}b^{*}\sO_{\mP(U^{5})}(1)\cdot\delta=0$,
$\sigma(\delta)$ is the $b$-exceptional curve over a point $\mathsf{s\in\gamma_{\mathtt{C}}}$.
Let $\mathsf{P}_{\mathsf{s}}:=b^{-1}(\mathsf{s})\simeq\mP^{2}$. Note
that the $U^{5}$-part of the coordinates of the point $\psi(\delta)$
is parallel to the coordinates of $\mathsf{s}$ by the description
of $\psi$-fibers as in the end of the proof of (1). Therefore, for
any $\psi$-exceptional curve $\delta'$ with $\psi(\delta')=\psi(\delta)$,
we have $b\circ\sigma(\delta)=b\circ\sigma(\delta')$. This implies
that, over $\psi(\mP(\sE^{\perp}|_{\mathsf{P}_{\mathsf{s}}}))$, the
$\psi$-fibers coincide with the corresponding fibers of $\mP(\sE^{\perp}|_{\mathsf{P}_{\mathsf{s}}})\to\psi(\mP(\sE^{\perp}|_{\mathsf{P}_{\mathsf{s}}}))$.
Since $\mathsf{P}_{\mathsf{s}}$ is $\mP^{2}$ and is isomorphically
mapped to a plane in $A_{\mathtt{C}}$, we have $\sE^{\perp}|_{\mathsf{P_{\mathsf{s}}}}=\Omega_{\mP^{2}}^{1}(1)\oplus\sO_{\mP^{2}}^{\oplus6}$,
and then $\psi$ induces the surjective map $\mP(\sE^{\perp}|_{\mathsf{P}_{\mathsf{s}}})\to\mP^{8}$.
It is easy to see that this is a $\mP^{1}$-bundle outside the image
$\mathsf{R_{\mathsf{s}}}\simeq\mP^{5}$ of $\mP(0\oplus\sO_{\mP^{2}}^{\oplus6})$,
over which the fibers of $\mP(\sE^{\perp}|_{\mathsf{P}_{\mathsf{s}}})\to\mP^{8}$
are $\mP^{2}$. Let $\mathsf{R}:=\cup_{\mathsf{s\in\gamma_{\mathtt{C}}}}\mathsf{R}_{\mathsf{s}}$,
which is a 6-dimensional variety. We have seen that $\delta$ coincides
with the $\psi$-fiber if $\psi(\delta)\in\psi(F_{\mP(\sE^{\perp})})\setminus(\psi(E_{\mP(\sE^{\perp})})\cup\mathsf{R}).$
Note that by a standard computation, we have $-K_{\mP(\sE^{\perp})}=8H_{\mP(\sE^{\perp})}+\sigma^{*}a^{*}\sO_{A_{\mathtt{C}}}(1)$.
Thus, by $\sigma^{*}a^{*}\sO_{A_{\mathtt{C}}}(1)\cdot\delta>0,$ we
have $-K_{\mP(\sE^{\perp})}\cdot\delta>0$. Therefore, $\overline{\Sigma}^{*}$
is smooth and $\psi$ is the blow-up along $\psi(F_{\mP(\sE^{\perp})})$
outside $\psi(E_{\mP(\sE^{\perp})})\cup\mathsf{R}$ by the proof of
\cite[Thm.2.3]{An}. By the description of the singular locus of $\overline{\Sigma}^{*}$
as in (2), we see that $S_{F}$ defined as in the statement of (2)
must be contained in $\mathsf{R}.$ Since $S_{F}$ and $\mathsf{R}$
are 6-dimensional and $\mathsf{R}$ is irreducible, we have $\mathsf{R}=S_{F}$.
Therefore, we obtain the first assertion.

The second and third assertions follow since the restriction of $\psi$
over $\mP((\wedge^{2}V^{3})^{*}\oplus(U^{5})^{*}\oplus0)$ is the
natural morphism $\mP(a^{*}(\Omega_{\mP^{7}}^{1}(1))\oplus0)\to\mP((\wedge^{2}V^{3})^{*}\oplus(U^{5})^{*}\oplus0)$
which is nothing but the universal family of the total transforms
of hyperplane sections of $A_{\mathtt{C}}.$
\end{proof}

\subsection{Cubic 3-fold and 4-fold}
\begin{cor}[Cubic 3-fold]
\label{cor:Cubic 3-fold}  Any smooth cubic $3$-fold is a linear
section of the cubic $\overline{\Sigma}^{*}$.
\end{cor}

\begin{proof}
We take $\Lambda$ as in Theorem \ref{thm:main2}. Then $\mP(\sE^{\perp})_{\Lambda}\to b\circ\sigma(\mP(\sE^{\perp})_{\Lambda})$
can be identified with the blow-up $Y'\to B_{3}$ along the twisted
cubic curve $C$. Note that $Y'$ has exactly two non-trivial contractions,
one of which is $Y'\to B_{3}$ and another is the anti-canonical morphism
$Y'\to W$. Since $\overline{\Sigma}^{*}\cap\mP(\Lambda)$ is a cubic
3-fold by Proposition \ref{prop:Dual Cubic}, $\mP(\sE^{\perp})_{\Lambda}\to\overline{\Sigma}^{*}\cap\mP(\Lambda)$
must coincide with $\mP(\sE^{\perp})_{\Lambda}\to b\circ\sigma(\mP(\sE^{\perp})_{\Lambda})$.
Now the assertion follows since any cubic 3-fold appears as $X'$
by \cite[II, Proof of Thm.0.10 (B) and (C)]{Tak1}. 
\end{proof}
\begin{cor}[Cubic 4-fold]
\label{corcubic 4-fold} Let $\Lambda\subset(\wedge^{2}V^{3})^{*}\oplus(U^{5})^{*}\oplus U^{5}$
be a general linear subspace of dimension $6$. The following assertions
hold: 

\vspace{3pt}

\noindent $(1)$ $\overline{\Sigma}^{*}\cap\mP(\Lambda)$ is a cubic
$4$-fold with one ordinary double point $\mathsf{v\in\mP}((\wedge^{2}V^{3})^{*}\oplus(U^{5})^{*}\oplus0)$.

\vspace{3pt}

\noindent $(2)$ Outside of $\mathsf{v}$, the induced morphism $\psi|_{\mP(\sE^{\perp})_{\Lambda}}\colon\mP(\sE^{\perp})_{\Lambda}\to\overline{\Sigma}^{*}\cap\mP(\Lambda)$
is the blow-up along $\psi(F_{\mP(\sE^{\perp})})\cap\mP(\Lambda)$
which is a twisted cubic cone with $\mathsf{v}$ as the vertex. The
$\psi|_{\mP(\sE^{\perp})_{\Lambda}}$-fiber over $\mathsf{\mathsf{v}}$
is isomorphic to the smooth 3-fold obtained by blowing up $B_{5}$
along a line.

\vspace{3pt}

\noindent $(3)$ Let $T_{\Lambda}:=\pi|_{\mP(\sE)_{\Lambda}}(\mP(\sE)_{\Lambda})$
. The morphism $\pi|_{\mP(\sE)_{\Lambda}}\colon\mP(\sE)_{\Lambda}\to T_{\Lambda}$
is an isomorphism and $T_{\Lambda}$ is a smooth $K3$ surface which
is isomorphic to a complete intersection in $A_{\mathtt{C}}$ of a
quadric containing $\Pi$ and a hyperplane.

\vspace{3pt}

\noindent $(4)$ The morphism $\sigma|_{\mP(\sE^{\perp})_{\Lambda}}\colon\mP(\sE^{\perp})_{\Lambda}\to\widehat{A}_{\mathtt{C}}$
is the blow-up of $\widehat{A}_{\mathtt{C}}$ along $T_{\Lambda}$. 
\end{cor}

\begin{proof}
The assertion (1) follows from Proposition \ref{prop:Dual Cubic}
(1) and (2), and the assertion (2) follows from Proposition \ref{prop:Dual Cubic}
(1) and (4).\vspace{3pt}

\noindent (3). Let $\Lambda'\subset(\wedge^{2}V^{3})^{*}\oplus(U^{5})^{*}\oplus U^{5}$
be a general linear subspace of dimension 7 containing $\Lambda$.
By \cite[Cor.5.18]{Tak2}, the restriction $\mP(\sE)_{\Lambda'}\to\overline{\Sigma}\cap\mP((\Lambda')^{\perp})$
of $\mP(\sE)\to\overline{\Sigma}$ can be identified with $Y'\to W$.
Note that we have $H_{\widehat{\Sigma}}|_{\mP(\sE)_{\Lambda'}}=\pi^{*}a^{*}\sO_{A_{\mathtt{C}}}(1)|_{\mP(\sE)_{\Lambda'}}$
by \cite[(5.7) and Lem.5.16]{Tak2}. Therefore $\mP(\sE)_{\Lambda'}\to\overline{\Sigma}\cap\mP((\Lambda')^{\perp})$
can be identified with $\mP(\sE)_{\Lambda'}\to W$ where $W$ is regarded
as a quadric section of $A_{\mathtt{C}}$ containing $\Pi$. Moreover
we may consider $\mP(\sE)_{\Lambda}$ is a general member of $|\pi^{*}a^{*}\sO_{A_{\mathtt{C}}}(1)|_{\mP(\sE)_{\Lambda'}}|$
and hence the image $T'_{\Lambda}\subset A_{\mathtt{C}}$ of $\mP(\sE)_{\Lambda}$
on $A_{\mathtt{C}}$ is a complete intersection in $A_{\mathtt{C}}$
of a quadric containing $\Pi$ and a hyperplane. Since $T'_{\Lambda}$
is disjoint from exceptional curves of $Y'\to W$ by generality, we
see that $\mP(\sE)_{\Lambda}\to T{}_{\Lambda}\to T_{\Lambda}'$ is
an isomorphism and hence $T_{\Lambda}\simeq T'_{\Lambda}$ is a smooth
$K3$ surface.

\vspace{3pt}

\noindent (4). Let $\mathsf{p}$ be a point of $\widehat{A}_{\mathtt{C}}$.
Considering the case that $l=\dim\Lambda=6$ and $r=5$ in the setting
of Lemma \ref{lem:HoTak}, we have $\dim(\sE_{\mathsf{p}}\cap\Lambda^{\perp})+1=\dim(\sE_{\mathsf{p}}^{\perp}\cap\Lambda).$
This implies that $\sigma|_{\mP(\sE^{\perp})_{\Lambda}}$ has nontrivial
fibers only over $T_{\Lambda}$ and they are isomorphic to $\mP^{1}$
by (3). Since $-K_{\mP(\sE^{\perp})_{\Lambda}}=(H_{\mP(\sE^{\perp})}+\sigma^{*}a^{*}\sO_{A_{\mathtt{C}}}(1))|_{\mP(\sE^{\perp})_{\Lambda}}$,
this is relatively ample over $\widehat{A}_{\mathtt{C}}.$ We see
that the relative Picard number of the morphism $\sigma|_{\mP(\sE^{\perp})_{\Lambda}}$
is one by the description of the fibers. Therefore, by \cite[Thm.2.3]{An},
$\sigma|_{\mP(\sE^{\perp})_{\Lambda}}$ is the blow-up of $\widehat{A}_{\mathtt{C}}$
along $T_{\Lambda}$.
\end{proof}
We immediately see that a cubic 4-fold $R$ with a double point $\mathsf{t}$
is rational projecting it from $\mathsf{t}$, and if $R$ is general,
then the blow-up of $R$ at $\mathsf{t}$ is equal to the blow-up
of $\mP^{4}$ along a smooth $K3$ surface which is a complete intersection
of a quadric and a cubic in $\mP^{3}$ (see \cite[Sect.5]{Ku5} for
further discussions). We have seen in Corollary \ref{corcubic 4-fold}
that if $R$ is a special one containing a cone over a twisted cubic,
then $R$ has another birational model which can be realized as the
blow-up along a $K3$ surface. 


\section{\textbf{$\mQ$-Fano 3-fold of genus 8}}

\subsection{Descriptions of $\mP(\sE^{\perp})$ \label{subsec:DescriptionsPEperpG8}}
\begin{prop}
\label{prop:Genus 8 dual} The following assertions hold:

\vspace{3pt}

\noindent $(1)$ The morphism $\psi$ is surjective and decomposes
as 
\[
\mP(\sE^{\perp})\stackrel{\psi_{1}}{\longrightarrow}\overline{\mP}\stackrel{\psi_{2}}{\longrightarrow}\mP(U^{3}\oplus(U^{3})^{*}\oplus\mathrm{S}^{-1,0,1}U^{3})
\]
where the morphism $\psi_{1}$ is birational and crepant, and $\psi_{2}$
is a finite morphism of degree $2$ branched along a sextic hypersurface
$\sB$. 

\vspace{3pt} 

\noindent $(2)$ The $\psi$-image of the $\psi_{1}$-exceptional
locus is the singular locus of $\sB$. 
\end{prop}

\begin{proof}
It is well-known that
\begin{align*}
c_{t}(T_{\mP(U^{7})}(-1)|_{B_{5}}) & =1+c_{1}(\sO_{B_{5}}(1))t+c_{1}(\sO_{B_{5}}(1))^{2}t^{2}+c_{1}(\sO_{B_{5}}(1))^{3}t^{3}\\
 & =1+c_{1}(\sO_{B_{5}}(1))t+(5l)t^{2}+5t^{3},
\end{align*}
where $l$ is the class of a line in $B_{5}.$ By \cite[Ex.3.2]{AC}
for example, we have $c_{t}(\sU)=1-c_{1}(\sO_{B_{5}}(1))t+(2l)t^{2},$
and hence the restriction of the universal exact sequence $0\to\sU|_{B_{5}}\to V'\otimes\sO_{B_{5}}\to\sQ|_{B_{5}}\to0$
gives 
\[
c_{t}(\sQ)=1+c_{1}(\sO_{B_{5}}(1))t+(3l)t^{2}+t^{3}.
\]
Therefore, we obtain, by a standard computation, 
\[
c_{t}((\sE^{\perp})^{*})=1+c_{1}(\sO_{B_{5}}(2))t+(13l)t^{2}+14t^{3}.
\]
Finally, we obtain 
\[
H_{\mP(\sE^{\perp})}^{11}=s_{3}((\sE^{\perp})^{*})=c_{1}(\sO_{B_{5}}(2))^{3}-2(c_{1}(\sO_{B_{5}}(2))\cdot(13l)+14=2.
\]
Therefore, since $\dim\mP(\sE^{\perp})=\dim\mP((V')^{*}\oplus(U^{7})^{*})=11$
and $-K_{\mP(\sE^{\perp})}=9H_{\mP(\sE^{\perp})},$ the assertion
(1) follows (the decomposition of $\psi$ is nothing but the Stein
factorization). Since $\psi_{1}$ is crepant, the singular locus of
$\overline{\mP}$ coincides with the $\psi_{1}$-image of the $\psi_{1}$-exceptional
locus. Thus the assertion (2) follows from a standard property of
the branched locus of a finite double cover.
\end{proof}

\subsection{Curve of genus $2$}
\begin{cor}
\label{cor:genus 2} Let $\Lambda$ be a $2$-dimensional subspace
of $U^{3}\oplus(U^{3})^{*}\oplus\mathrm{S}^{-1,0,1}U^{3}$. If $\mP(\Lambda)\cap\sB$
consists of exactly $6$ points, then $\mP(\sE^{\perp})_{\Lambda}\to\mP^{1}=\mP(\Lambda)$
is a finite double cover branched along $6$ points. 
\end{cor}

\begin{proof}
If $\dim\Lambda$=2 and $\mP(\Lambda)\cap\sB$ consist of exactly
$6$ points, $\mP(\Lambda)$ is disjoint from the singular locus of
$\sB$. Therefore the assertion follows from Proposition\ref{prop:Genus 8 dual}.
\end{proof}
%


As a general property of a curve of genus 2, we have the following:
\begin{prop}
For any smooth curve $C$ of genus $2$ and any divisor $\delta$
on $C$ of degree $7$, there exists a prime $\mQ$-Fano $3$-fold
$X$ of genus $8$ such that $f'\colon Y'\to X'=B_{5}$ is the blow-up
along a curve isomorphic to $C$ and $\sO_{B_{5}}(1)$ restricts to
$\delta$. 
\end{prop}

\begin{proof}
Let $\varepsilon:=\delta-K_{C}$, which is a divisor of degree 5 and
hence is very ample. We consider that $C$ is embedded in $\mP^{3}$
by $|\varepsilon|$. We choose a projection of $\mP^{3}\dashrightarrow\mP^{2}$
from a point outside of $C$ such that the image $\overline{C}$ of
$C$ is a quintic 4-nodal plane curve. Since a line cannot pass through
3 nodes of $\overline{C}$, the 4 nodes of $\overline{C}$ is in a
general position. Let $T\to\mP^{2}$ be the blow-up at the 4 nodes
of $\overline{C}$, and $m$ and $e_{i}\,(1\leq i\leq4)$ are the
total transform of a line on $\mP^{2}$ and the exceptional curves
of the blow-up respectively. Note that $T$ is a smooth quintic del
Pezzo surface and hence we consider $T$ is a hyperplane section of
$B_{5}$. The strict transform of $\overline{C}$ is smooth, hence
we denote it by $C$. Then $C\sim5m-2\sum_{i=1}^{4}e_{i}$ and it
holds that  $\varepsilon=m|_{C}$ and $K_{C}=(2m-\sum_{i=1}^{4}e_{i})|_{C}$.
Hence we have $\delta=K_{C}+\varepsilon=(3m-\sum_{i=1}^{4}e_{i})|_{C}=-K_{T}|_{C}=\sO_{B_{5}}(1)$.
Now, by the proof of \cite[Part II, Thm.0.10 (B)]{Tak1}, we see that
the blow-up of $B_{5}$ along $C$ is a part of the Sarkisov link
(\ref{eq:Sarkisov}) for a prime $\mQ$-Fano $3$-fold $X$of genus
$8$.
\end{proof}

\section{\textbf{Trinity}}

Finally, with some compensations, we sum up three types of appearances
of the curve $C$ appearing in the basic diagram except in the case
of genus 6 and $\mathtt{C}$-type:
\begin{thm}
\label{thm:trinity} We fix one of the $5$ classes of $\mQ$-Fano
$3$-fold except the class of genus $6$ and $\mathtt{C}$-type, and
let $\sE$ be the vector bundle as in Table 1 for the class. The following
assertions are equivalent for a smooth curve $C$:

\vspace{3pt}

\noindent $(1)$ For a linear subspace of $V_{\sE}^{*}$ of dimension
$r-1$, $\mP(\sE)_{\Lambda}$ appears as $Y'$ or $Z'$ in the basic
diagram and $\mP(\sE)_{\Lambda}\to S$ is the blow-up along a curve
isomorphic to $C$.

\vspace{3pt}

\noindent $(2)$ $C\simeq\sigma(\mP(\sE^{\perp})_{\Lambda})\simeq\mP(\sE^{\perp})_{\Lambda}$
for a linear subspace of $V_{\sE}^{*}$ of dimension $r-1$.

\vspace{3pt}

\noindent $(3)$ For a linear subspace of $V_{\sE}^{*}$ of dimension
$r-1$, $C$ is the double cover of $\mP(\Lambda)$ branched along
$\sB\cap\mP(\Lambda)$ in the genus 8 case, or $C\simeq\overline{\Sigma}^{*}\cap\mP(\Lambda)$
in the other cases.
\end{thm}

\begin{proof}
\vspace{3pt}

\noindent $(1)\Rightarrow(2)$ This is proved in Theorem \ref{thm:main2}
(2).

\vspace{3pt}

\noindent $(2)\Rightarrow(1)$ Actually we only need the assumption
that $C\simeq\mP(\sE^{\perp})_{\Lambda}$. Since $\mP(\sE^{\perp})_{\Lambda}$
has the expected dimension, $\mP(\sE)_{\Lambda}$ has also the expected
dimension by Lemma \ref{lem:HoTak} with $r-l=1$. Since $\mP(\sE^{\perp})_{\Lambda}$
is smooth, so is $\mP(\sE)_{\Lambda}$ by \cite[Thm.7.12]{Ku4}. By
Lemma \ref{lem:HoTak} with $r-l=1$ again, nontrivial fibers of $\mP(\sE)_{\Lambda}\to S$
are $\mP^{1}$'s over $\sigma(\mP(\sE^{\perp})_{\Lambda}).$ Since
$-K_{\mP(\sE)_{\Lambda}}=H_{\mP(\sE)}|_{\mP(\sE)_{\Lambda}}$, we
see that $-K_{\mP(\sE)_{\Lambda}}$ is relatively ample over $S$.
Therefore, by \cite[Thm.2.3]{An}, the morphism $\mP(\sE)_{\Lambda}\to S$
is the blow-up along $\sigma(\mP(\sE^{\perp})_{\Lambda}).$ 

Since we see that $-K_{\mP(\sE)_{\Lambda}}=H_{\mP(\sE)}|_{\mP(\sE)_{\Lambda}}$
and $(-K_{\mP(\sE)_{\Lambda}})^{3}=2g(X)-2>0,$ $\mP(\sE)_{\Lambda}$
is a smooth weak Fano 3-fold. Restricting the diagram \cite[(3.2), (5.10), or (6.7)]{Tak2},
we can verify the assertion.

\vspace{3pt}

\noindent $(2)\Rightarrow(3)$ By the assumption (2), $C$ is the
normalization of the double cover of $\mP(\Lambda)$ branched along
$\mP(\Lambda)\cap\sB$ in the genus 8 case by Proposition \ref{prop:Genus 8 dual},
or the normalization of $\overline{\Sigma}^{*}\cap\mP(\Lambda)$ in
each of the other cases by Proposition \ref{prop:DualG4}, \ref{prop:G5 Dual}
or \cite[Prop.4.12]{Tak2} with explanations as in the section \ref{sec:genus 6Q}.
According to case-by-case check, the arithmetic genus of the double
cover of $\mP(\Lambda)$ branched along $\mP(\Lambda)\cap\sB$ in
the genus 8 case (resp. the 1-dimensional linear section $\overline{\Sigma}^{*}\cap\mP(\Lambda)$
in each of the other cases) is the same as the genus of $C$. Therefore
$C\simeq\overline{\Sigma}^{*}\cap\mP(\Lambda)$.

\vspace{3pt}

\noindent $(3)\Rightarrow(1)$ Since $\overline{\Sigma}^{*}\cap\mP(\Lambda)$
is smooth and is a 1-dimensional linear section of $\overline{\Sigma}^{*}$,
$\mP(\Lambda)$ is disjoint from $\Sing\overline{\Sigma}^{*}.$ Therefore
$\mP(\sE^{\perp})_{\Lambda}\simeq\overline{\Sigma}^{*}\cap\mP(\Lambda)\simeq C$.
Note that $(2)\Rightarrow(1)$ holds by the weaker assumption that
$\mP(\sE^{\perp})_{\Lambda}\simeq C$ as we have remarked in the proof.
Therefore $(3)\Rightarrow(1)$ follows.
\end{proof}

\end{document}